\documentclass{article}

\usepackage[utf8]{inputenc}
\usepackage{graphicx}
\usepackage{amsmath}
\usepackage{amsfonts}
\usepackage{amssymb}
\usepackage{amsthm}
\usepackage{a4wide}
\usepackage{stmaryrd}
\usepackage{siunitx}

\newcommand{\Xv}{{X}}

\newcommand{\Tv}{T}
\newcommand{\Nv}{N}

\newcommand{\uv}{u}
\newcommand{\nv}{n}
\newcommand{\tv}{t}
\newcommand{\fv}{f}

\newcommand{\xv}{{x}}

\newcommand{\vv}{v}
\newcommand{\alphav}{\alpha}
\newcommand{\sigmav}{{\sigma}}

\newcommand{\sv}{s}
\newcommand{\Sv}{S}
\newcommand{\gv}{g}

\newcommand{\It}{\mathbf{I}}
\newcommand{\Pt}{\mathbf{P}}
\newcommand{\Ft}{\mathbf{F}}
\newcommand{\Gt}{\mathbf{G}}
\newcommand{\Ct}{\mathbf{C}}
\newcommand{\Et}{\mathbf{E}}
\newcommand{\Mt}{\mathbb{D}}

\newcommand{\sigmat}{\boldsymbol{\sigma}}

\newcommand{\epst}{\boldsymbol{\varepsilon}}
\newcommand{\Sigmat}{\boldsymbol{\Sigma}}

\newcommand{\opdiv}{\operatorname{div}}
\newcommand{\opcurl}{\operatorname{curl}}
\newcommand{\skw}{\operatorname{skw}}
\newcommand{\sym}{\operatorname{sym}}

\newcommand{\T}{\mathcal{T}}
\newcommand{\E}{\mathcal{E}}
\newcommand{\Lag}{\mathcal{L}}

\newcommand{\imt}{{0}}
\newcommand{\upd}{{\Delta}}

\usepackage{chngcntr}
\counterwithin{figure}{section}
\counterwithin{table}{section}
\counterwithin{equation}{section}

\title{Three-field mixed finite element methods for nonlinear elasticity}

\author{\normalsize{
		Michael Neunteufel \thanks{Corresponding author: Michael Neunteufel, Institute for Analysis and Scientific Computing, Technische Universit\"at Wien, Wiedner Hauptstr. 8-10 , 1040 Wien, Austria, email: michael.neunteufel@tuwien.ac.at}
		\quad , \quad
		Astrid S. Pechstein \thanks{Astrid S. Pechstein, Institute of Technical Mechanics, Johannes Kepler University Linz, Altenbergerstr. 69, 4040 Linz, Austria, email: astrid.pechstein@jku.at}
		\quad
	}
	\normalsize{and \quad
		Joachim Sch\"oberl \thanks{Joachim Sch\"oberl, Institute for Analysis and Scientific Computing, Technische Universit\"at Wien, Wiedner Hauptstr. 8-10 , 1040 Wien, Austria, email: joachim.schoeberl@tuwien.ac.at}
	}
}
\date{}

\begin{document}
		
\maketitle
\begin{abstract}
In this paper, we extend the tangential-displacement normal-normal-stress
continuous (TDNNS) method from \cite{PS11} to nonlinear elasticity. By means of the Hu--Washizu
principle, the distibutional derivatives of the displacement vector
are lifted to a regular strain tensor. We introduce three different
methods, where either the deformation gradient, the Cauchy--Green
strain tensor, or both of them are used as independent variables.
Within the linear sub-problems, all stress and strain variables can be
locally eliminated leading to an equation system in displacement
variables, only.
The good performance and accuracy of the presented methods are
demonstrated by means of several numerical examples (available via www.gitlab.com/mneunteufel/nonlinear\_elasticity).\\

\end{abstract}

\section*{Introduction} 
\label{sec:intro}
The construction of discretization methods in elasticity is an active
field of research. 
Mixed finite element methods entail a large variety of beneficial attributes 
compared to standard finite elements, which suffer from volume locking in the 
nearly incompressible regime and shear locking for flat elements.
As an example of a mixed method, the Hellinger--Reissner formulation introduces
the stress as a second, independent variable. 
For linear problems, there is a well established mathematical framework for
analyzing stability and robustness properties,
see the monograph by Boffi, Brezzi and Fortin \cite{BBF13}. 
The application of Hellinger--Reissner mixed methods to nonlinear problems is limited.
Another classical mixed approach, the so-called Hu--Washizu or three-field approach, dates back to \cite{W68}.
Here, the strain is used as a third field, which allows to treat a large
class of nonlinear problems. 

A basic principle of mixed formulations is integration by parts,
which moves derivatives from the displacement to the stress
field. Thus, higher continuity of the approximation space for the
stress is required. The construction of symmetric matrix-valued finite
elements with continuous normal components is non-trivial.
Simplicial elements for the mixed Hellinger--Reissner formulation in two and three space dimensions have been constructed 
by Arnold and Winther \cite{AW02} or Arnold, Awanou and Winther \cite{AAW08}, respectively.
As, especially in three dimensions, these elements are rather costly, different efficient variants have been introduced as elements with weak symmetry (e.g. \cite{ABD84,Sten88,AFW07}).
In \cite{PS11,PS12,PS18}, the authors introduced and analyzed a Hellinger--Reissner formulation employing normal-normal continuous stress elements.
This \emph{tangential-displacement normal-normal-stress} (TDNNS) method for linear elasticity has been shown to be free from shear locking when discretizing thin structures by flat elements, and also robust when approaching the incompressible limit.

In Hu--Washizu formulations, three independent fields are introduced -- in the linear case, these three fields are comprised of displacement, strain, and stress, see also \cite{Braess07}. Different three-field formulations using an additional \emph{enhanced strain} have been proposed in the literature: we cite the works by Simo and Rifai \cite{SR90} and Kasper and Taylor \cite{KT00a}. Therein, the strain is assumed to decompose additively into a compatible part, associated with the displacement field, and an enhanced part. These methods were proven to be stable and efficient by various patch tests and benchmark problems. Reddy and Simo \cite{RS95} also provided rigorous stability and convergence proofs for their elements.

The enhanced strain method from \cite{SR90} has been extended to nonlinear elasticity by Simo, Armero and Taylor \cite{SA92,SAT93}. 
Earlier, a mixed formulation based on the Hu--Washizu principle including finite-deformation elasto-plasticity was proposed by Simo, Taylor and Pister \cite{STP85}.
Kasper and Taylor \cite{KT00b} proposed a non-linear mixed-enhanced method with independent deformation gradient and first Piola--Kirchhoff stress.  More recently, different groups introduced methods based on polyconvex strain energy-functions, where the deformation gradient, cofactor matrix and determinant are treated separately: we cite Pfefferkorn and Betsch \cite{PB20}, Schr\"oder, Wriggers and Balzani \cite{SWB11}, and Bonet, Gil and Ortigosa \cite{BGO15,BGO16}.

Reese, Wriggers and Reddy \cite{RWR98} introduced a reduced-integration stabilized brick element for finite elasticity, where the stabilization is based on the enhanced strain method. Later, Reese \cite{Reese07} proposed a brick solid-shell element based on reduced integration and hourglass stabilizations. More recently, a group around Reese proposed  lowest-order locking-free hybrid discontinuous Galerkin elements for large deformations \cite{WBAR17,BKWWRWW18}.  The equivalence of this approach to hourglass stabilization and reduced integration was discussed in \cite{RBW17}. 

Angoshtari, Shojaei and Yavari \cite{ASY17} introduced the compatible-strain mixed finite element method (CSMFEMs) for two-dimensional compressible large deformation problems, which also belongs to the class of three-field formulations. Shojaei and Yavari later generalized this method to three dimensions and incompressible elasticity in \cite{SY18,SY19}. These elements are based on a 
Hilbert complex of nonlinear elasticity \cite{AY16}: the displacement is chosen to be discretized with (conforming) nodal elements, its gradient by N{\'e}d{\'e}lec and the first Piola--Kirchhoff stress tensor by Raviart--Thomas elements.

Recently, virtual element methods (VEM) for nonlinear (in-)elastic problems in the compressible and incompressible regime have been proposed \cite{BLM15,ABLS17,CBP17,WRRH17}.

The linear TDNNS method introduced and analysed in \cite{PS11, PS18} does not suffer from shear locking or large aspect ratios \cite{PS12}. Furthermore, volume locking in the nearly incompressible case can be avoided by adding a consistent stabilization term \cite{Sinw09}. Therein, the displacement field is discretized by N{\'e}d{\'e}lec elements resulting in less regularity requirements as for nodal elements. This, however, prevents the use of nonlinear material laws. Recently an Updated Lagrangian (UL) mixed scheme has been proposed for non-linear elasticity \cite{P19}. This mixed method is a two-field approach, where at convergence the independent fields are displacement and Cauchy stress.

In this work, we choose a different approach: we use the Hu--Washizu three-field principle. In the following, two different variants of three-field methods are deployed extending the (linear) TDNNS method. 
In the first variant, deformation gradient and first Piola-Kirchhoff stress are introduced as independent quantities. The independent deformation gradient may be interpreted as a lifting of the distributional gradient of the (discontinuous) displacement field. This method proved very robust in our numerical examples.
In the second variant, Cauchy--Green strain and second Piola--Kirchhoff stress comprise the additional independent variables. In computations, rapid convergence of nonlinear iterations was observed.
In a third approach, five independent fields are used: displacement, deformation gradient, Cauchy--Green strain, and first and second Piola--Kirchhoff stress. Mathematically, this five-field formulation  corresponds to a lifting followed by a projection, trying to combine the advantages of the two three-field methods. The additional fields are discretized by a discontinuous version of so-called Regge elements; originally derived as geometric discretization of the Einstein field equation by Regge \cite{Regge61}, further developed in \cite{cheg86}, and given a finite element context through Christiansen \cite{christiansen11}; enabling static condensation techniques. Therefore the final system involves only displacement and stress degrees of freedom making the methods computationally efficient.

This paper is organized as follows: in Section~\ref{sec:largedeformation} notation of elasticity is introduced and Section~\ref{sec:smalldefTDNNS} is devoted to summarize the linear TDNNS method. In Section~\ref{sec:nonlinTDNNS} the three extensions of the TDNNS method to the large deformation regime are presented and an Updated Lagrangian scheme for these approaches is derived in Section~\ref{sec:updated_lagrangian}. Finally, Section~\ref{sec:numerics} describes the used finite element spaces in more detail and several numerical examples are presented, showing the excellent performance of the methods.

\section{Large deformation elasticity} \label{sec:largedeformation}
Let $\Omega \subset \mathbb R^d$ with $d=2,3$ denote the body of interest in (stress-free) reference configuration. 
We use reference coordinates $\Xv \in \Omega$ to describe the position of a material point. Under deformation, 
material points $\Xv$ are mapped to spatial positions $\xv(\Xv)$. The associated deformation field shall be denoted by
$\uv$ and is defined as
\begin{align}
&\uv: \Omega \to \mathbb R^d, & \uv(\Xv) &= \xv(\Xv) - \Xv.
\end{align}
In this work, all spatial derivatives are assumed to be with respect to $\Xv$ if not indicated otherwise. 
We use $\nabla = \nabla_{\Xv}$, $\opcurl = \opcurl_{\Xv}$ and $\opdiv = \opdiv_{\Xv}$ defined in the standard ways.
When applied to a second order tensor, the divergence operator is to be understood row-wise.

The deformation gradient $\Ft$ is defined as
\begin{align}
\Ft := \It + \nabla \uv.
\end{align}
We use $J := \det \Ft$ for the Jacobi determinant. Further, we introduce the (right) Cauchy--Green strain tensor $\Ct$ and the Green strain tensor $\Et$ by
\begin{align}
\Ct := \Ft^\top\Ft, && \Et:=\frac{1}{2}\left(\Ct-\It\right). \label{eq:cg_g_tensor}
\end{align} 

Throughout this work, we
assume to be in a static regime without inertial forces. 
On a non-vanishing part $\Gamma_D$ of the boundary,
the displacement shall be prescribed such that $\uv = \uv_D$ on $\Gamma_D$. On the remaining part $\Gamma_N = \partial \Omega \backslash \Gamma_D$, surface tractions $\gv$ are given. Moreover, the body of interest shall be subjected to an external body force $\fv$. In an energy-based formulation, we need the work of external forces, which sums up to
\begin{align}
W_{ext} &= \int_\Omega \fv\cdot\uv\,d\Xv + \int_{\Gamma_N}\gv\cdot\uv\,d\Sv.
\end{align}

Let $\Psi(\cdot)$ be a hyperelastic potential, where we use the same symbol independently of which deformation measure $\Ft$, $\Ct$, or $\Et$ is used.

The (non-symmetric) first Piola--Kirchhoff stress is defined as the derivative
\begin{align}
\Pt := \frac{\partial \Psi}{\partial \Ft}.\label{eq:pk1}
\end{align}
The symmetric second Piola--Kirchhoff stress tensor is related to $\Psi$ and $\Pt$ via
\begin{align}
\Sigmat &:= 2\frac{\partial \Psi}{\partial \Ct}, & \text{and }&&\Ft \Sigmat = \Pt .\label{eq:pk2}
\end{align}
For the sake of completeness, we introduce the (symmetric) Cauchy stress tensor $\sigmat$, which is obtained transforming the first Piola--Kirchhoff stress to spatial configuration by the Piola transform,
\begin{align}
\sigmat = \frac{1}{J} \Pt \Ft^\top.
\end{align}

The displacement field $\uv$ satisfies the minimization problem
\begin{align}
\int_{\Omega}\Psi(\Ct)\, d\Xv- W_{ext} &\to \min_{\uv \in V}, \label{eq:min_prob_nlinel}
\end{align}
where $V$ denotes the set of all admissible displacement fields. We do not give a precise mathematical definition of the set $V$ here, but only note that the condition ``$\uv = \uv_D$ on $\Gamma_D$'' is essential. Conditions on continuity and/or (weak) differentiability of $\uv \in V$ are detailed for the different finite element methods in Section~\ref{sec:smalldefTDNNS} and Section~\ref{sec:nonlinTDNNS}.

Computing the first variation of \eqref{eq:min_prob_nlinel} in direction $\delta \uv$ yields the following variational equation, well known as \emph{principle of virtual works}: Find $\uv\in V$ such that
\begin{align}
\int_{\Omega}\Pt:\nabla\delta\uv\,d\Xv = \int_{\Omega}\fv\cdot\delta\uv\,d\Xv + \int_{\Gamma_N}\gv\cdot\delta\uv\,d\Sv\quad\forall\delta \uv,\label{eq:var_prob_nlinel}
\end{align}
where $\delta \uv\in V_0$ lives in the set of admissible virtual displacements $V_0$ with $\delta \uv = 0$ on $\Gamma_D$.
Integration by parts leads to the balance equation, which reads in strong form
\begin{align}
-\opdiv(\Pt) = \fv, && \Pt \Nv =\gv \text{ on }\Gamma_N, && \uv=\uv_D \text{ on }\Gamma_D.\label{eq:balance}
\end{align}
Above, $\Nv$ denotes the outer normal vector to $\partial \Omega$ in reference configuration.

\paragraph{A linear small-deformation model} 

As it will be needed  for theoretical considerations in the sequel, we shortly present the according linearized
elasticity problem. Under the assumption of small deformations, material points $\Xv$ and spatial points $\xv$
can be identified, leading to the linearized strain tensor 
\begin{align}
\epst(\uv) := \frac12 (\nabla \uv + (\nabla \uv)^\top)
\end{align}
as a measure of the deformation. As the different stress tensors coincide in this case, we use $\sigmat$ for the stress in linearized elasticity. Assuming the potential $\Psi$ to be quadratic in $\epst$, the stress-strain relation is linear $\sigmat =\Mt\epst(\uv)$, represented by the fourth order stiffness tensor $\Mt$. 
The principle of virtual works transforms from \eqref{eq:var_prob_nlinel} to
\begin{align}
\int_{\Omega}\Mt\epst(\uv):\epst(\delta\uv)\,d\Xv = \int_{\Omega}\fv\cdot\delta\uv \,d\Xv + \int_{\Gamma_N}\gv\cdot\delta\uv\,d\Sv. \label{eq:var_prob_linel}
\end{align}

\paragraph{Basic finite element ingredients}
To introduce any finite element method, let $\T = \{T\}$ be a regular finite element mesh of $\Omega$. In the past, TDNNS elements for triangles, quadrilaterals, as well as tetrahedra, hexahedra and prismatic elements have been developed. The mesh $\T$ may be a hybrid one, containing several of the aforementioned element types. For each element $T$ we define its boundary $\partial T$ and the corresponding outer normal vector $\Nv$. 

We introduce the following conventions for denoting normal and tangential components of vector and tensor fields:
the normal component of a vector field $\vv$ is denoted by $v_{\Nv}:=\vv\cdot\Nv$, while the tangential component is given by $\vv_{\Tv}:= \vv - v_{\Nv} \Nv$. For a tensor field $\sigmat$, we introduce its normal vector $\sigmav_{\Nv} = \sigmat \Nv$.
This normal vector can again be split into a normal and tangential component, reading $\sigma_{\Nv\Nv} = (\sigmat \Nv)\cdot \Nv$ and $\sigmav_{\Nv\Tv} = \sigmat \Nv - \sigma_{\Nv\Nv} \Nv$, respectively.

\section{The linear TDNNS method} \label{sec:smalldefTDNNS}

The tangential-displacement and normal-normal-stress (TDNNS) method introduced for linear elasticity in \cite{PS11, PS12, PS18} uses mixed finite elements with tangential continuous displacement fields and normal-normal continuous stresses. These are also chosen as degrees of freedom (dofs). 

We introduce the linear TDNNS method shortly; for details on the finite element spaces we refer to Section~\ref{sec:numerics}.

Assume $\uv$ to be an admissible displacement finite element function, where by admissible we mean: $\uv$ is piecewise smooth on $\T$, its tangential component $\uv_{\Tv}$ is continuous across element interfaces, and the tangential component satisfies the essential boundary condition on $\Gamma_D$,  $\uv_{\Tv} =(\uv_D)_{\Tv}$.
As the normal components are not necessarily continuous, gaps in normal direction may open up at element interfaces,
and need to be treated accordingly. As the stress is treated as an independent variable, admissible stress fields
$\sigmat$ are introduced, where admissible means: $\sigmat$ is piecewise smooth symmetric on $\T$, its normal-normal component $\sigma_{\Nv\Nv}$ is continuous across element interfaces, and it satisfies the boundary condition on $\Gamma_N$: $\sigma_{\Nv\Nv} = \gv_{\Nv}$. The continuous normal component of the stress vector will control the displacement gaps.\\

Instead of the minimization problem \eqref{eq:min_prob_nlinel}, a saddle point problem for the Lagrangian is posed:
find admissible displacement and stress fields minimizing respectively maximizing the Lagrangian
\begin{align}
\mathcal L(\uv, \sigmat) := - \frac12  \int_\Omega \Mt^{-1} \sigmat:\sigmat\,d\Xv + \langle \epst(\uv) , \sigmat \rangle_\T - W_{ext}^{\text{TDNNS}} \to \min_{\uv \text{ adm.}} \max_{\sigmat \text{ adm.}}.\label{eq:var_prob_lin_tdnns}
\end{align}
Note that the Lagrangian corresponds to the mechanical enthalpy of the system. 
The duality pairing $\langle\cdot,\cdot\rangle_\T$ is defined in the sense of distributions as neither the divergence of the stress nor the gradient of the displacement field is a globally regular function \cite{PS11}
\begin{align}
\langle \epst(\uv), \sigmat \rangle_\T &:= \sum_{T\in\T} \Big( \int_T \sigmat:\nabla\uv\,d\Xv - \int_{\partial T} \sigmat_{\Nv\Nv}\uv_{\Nv}\,d\Sv\Big)\nonumber\\
&=
- \sum_{T\in\T} \Big( \int_T \text{div}(\sigmat)\cdot\uv\,d\Xv - \int_{\partial T} \sigmat_{\Nv\Tv}\cdot\uv_{\Tv}\,d\Sv \Big) 
= -\langle \text{div}(\sigmat),\uv\rangle_\T.
\label{eq:def_dual_pair2}
\end{align}
The work of external forces is adapted to read
\begin{align}
W^{\text{\text{TDNNS}}}_{ext} &= \int_\Omega \fv\cdot\uv\,d\Xv + \int_{\Gamma_N}\gv_{\Tv}\cdot\uv_{\Tv}\,d\Sv - \int_{\Gamma_D} (u_D)_{\Nv}\, \sigma_{\Nv\Nv}\,d\Sv. \label{eq:Wtdnns}
\end{align}

After discretization, the above saddle point problem leads to a linear system of equations with an indefinite system matrix.
However, a positive definite system matrix can be obtained using a \emph{hybridization} technique. Therefore, the normal-normal
continuity of the stresses is broken, such that $\sigma_{\Nv\Nv}$ may be discontinuous across interfaces, and does not
necessarily satisfy the stress boundary condition on $\Gamma_N$. A Lagrange multiplier $\alphav$ is introduced, which
enforces the lost continuity and boundary conditions. As, in the finite element scheme, both $\sigma_{\Nv\Nv}$ and
$\alphav$ are of the same polynomial order, the hybridized system is equivalent to the original one. 

To be  precise,
let $\alphav$ be a vector-valued function defined on the skeleton $\E$, whose direction is normal to the interfaces.
To be admissible, $\alphav$ shall moreover satisfy that $\alphav_{\Nv}= (\uv_D)_{\Nv}$ on the displacement 
boundary $\Gamma_D$. The hybridized TDNNS problem then reads

\begin{align}
&\mathcal L^{h}(\uv, \sigmat, \alphav) := \nonumber\\
&- \frac12  \int_\Omega \Mt^{-1} \sigmat:\sigmat\,d\Xv + \langle \epst(\uv) , \sigmat \rangle_\T + 
\sum_{T \in \T} \int_{\partial T} \sigma_{\Nv\Nv} \alpha_{\Nv}\, d\Sv - W_{ext}^{\text{TDNNSh}} \to 
\min_{\uv \text{ adm.}}  \min_{\alphav \text{ adm.}} \max_{\sigmat \text{ disc.}}, \label{eq:var_prob_lin_tdnns_hyb}
\end{align}
with the work of external forces
\begin{align}
W^{\text{TDNNSh}}_{ext} &= \int_\Omega \fv\cdot\uv\,d\Xv + \int_{\Gamma_N}( \gv_{\Tv}\cdot\uv_{\Tv} + g_{\Nv} \, \alpha_{\Nv})\,d\Sv.
\label{eq:WTDNNSh}
\end{align}

The additive term containing $\alphav$ indeed enforces normal-normal continuity of $\sigmat$, as can be seen by the following consideration:
if the surface integrals are re-ordered facet by facet, there are two contributions to each internal facet with normal 
vectors $\Nv$ of opposite direction. Thus, $\alpha_{\Nv} = \alphav \cdot \Nv$ changes its sign, while $\sigma_{\Nv\Nv}$
keeps its sign due to the quadratic occurrence of $\Nv$. Considering boundary facets correctly, using $\Nv_E$ as a unique facet normal, and denoting
the jump of the normal stress by $\llbracket\sigma_{\Nv_E\Nv_E}\rrbracket$, one observes
\begin{align}
\sum_{T \in \T} \int_{\partial T} \sigma_{\Nv\Nv} \alpha_{\Nv}\, d\Sv &= \sum_{E \in \E} \int_{E} \llbracket\sigma_{\Nv_E\Nv_E}\rrbracket \alpha_{\Nv_E}\, d\Sv,
\end{align}
where $\E$ denotes the set of all interfaces, i.e., edges in two dimensions and faces in 3D.

It shows that the Lagrange multiplier $\alphav$ has the physical meaning of the normal component of the displacement \cite{PS11}.
As the discontinuous stress $\sigmat$ in \eqref{eq:var_prob_lin_tdnns_hyb} does not have any coupling dofs, one can use static condensation to eliminate it at element level, reducing the number of total dofs drastically for the final system, and making it therefore symmetric and positive definite (spd) again.

\section{Nonlinear TDNNS} \label{sec:nonlinTDNNS}
For nonlinear and thus, in general not explicitly invertible material laws the TDNNS methods \eqref{eq:var_prob_lin_tdnns} and \eqref{eq:var_prob_lin_tdnns_hyb} can not be applied. The main difficulty is that the gradient of the displacement field $\uv$ is a distribution rather than a function. Therefore, multiplication might not be well defined, which is, however, crucial for handling nonlinear materials. In \cite{P19} recently an Updated Lagrangian scheme has been discussed to enable these sort of materials.

In this work, however, we will use the Hu--Washizu principle \cite{W68} introducing a new independent field. We will discuss three different approaches. In the first one we will lift $\nabla \uv +\It$ to the deformation gradient $\Ft$ as a new unknown, which will be a regular function again. The second ansatz uses that {for objective materials the energy potential $\Psi$ depends on the Cauchy--Green strain tensor} $\Ct$ - or equivalently on the Green strain tensor $\Et$ - taking it as an additional field. The third approach combines the first two by first introducing a lifting to $\Ft$ followed by a projection to $\Ct$, i.e., two additional fields are used.

In what follows, we will use the notation $\Ft(\uv) := \It + \nabla \uv$, $ \Ct(\uv):= \Ft(\uv)^\top\Ft(\uv)$ and $\Et(\uv):= 0.5(\Ct(\uv)-\It)$ to indicate the dependence on the displacement, whereas $\Ft$, $\Ct$, and $\Et$ denote independent fields. Further, without restriction of generality but to simplify the problem at hand, we will use homogeneous traction forces $\gv=0$ on $\Gamma_N$ and homogeneous Dirichlet data $\uv_D=0$ on $\Gamma_D$. 

\subsection{Lifting to $\Ft$}
\label{subsec:lifting_F}
Instead of solving \eqref{eq:min_prob_nlinel} one may define the following constrained minimization problem for piecewise smooth $\tilde \uv$ on $\Tv$, where $\tilde \uv$ is continuous and satisfies the displacement boundary condition on $\Gamma_D$,
\begin{align}
\int_{\Omega}\Psi(\Ft)\,d\Xv - W_{ext} &\to \min\limits_{\substack{\tilde \uv\text{ cont.}\\\Ft=\Ft(\tilde \uv)}}.\label{eq:min_prob_huwa_F}
\end{align}
The corresponding Lagrange functional reads
\begin{align}
\Lag(\uv, \Ft,\Pt):=\int_{\Omega}\Psi(\Ft)\,d\Xv - W_{ext} - \int_\Omega (\Ft-\It-\nabla\uv):\Pt\, d\Xv,\label{eq:lag_huwa_F}
\end{align}
with the first Piola--Kirchhoff stress tensor as Lagrange multiplier. For tangential-continuous displacement
functions, however, the last integral in \eqref{eq:lag_huwa_F} is not well-defined, as $\nabla \uv$ does not exist
in the sense of a square-integrable function. For normal-normal continuous (but non-symmetric) $\Pt$, the integral 
can be re-interpreted as a distribution. For $\uv$ piecewise smooth tangential-continuous, $\Pt$ piecewise smooth 
non-symmetric with $P_{\Nv\Nv}$ continuous, and $\Ft$ piecewise smooth discontinuous we obtain the
saddle point problem
\begin{align}
\Lag^{\Ft}(\uv, \Ft,\Pt) \to &\min_{\uv \text{ adm.}} \max_{\Pt \text{ adm.}} \min_{\Ft \text{ disc.}},\label{eq:var_prob_huwa_F}\\
\text{with } \Lag^{\Ft}(\uv, \Ft,\Pt) &:=\int_{\Omega}\Psi(\Ft)\,d\Xv - W_{ext}^{\text{TDNNS}} + \langle \nabla \uv, \Pt\rangle_\T - \int_\Omega (\Ft-\It):\Pt\, d\Xv\\
\text{and } \langle \nabla \uv, \Pt\rangle_\T &=
\sum_{T \in \T} \Big( \int_T \nabla \uv : \Pt\,d\Xv - \int_{\partial T} \uv_{\Nv} \, P_{\Nv\Nv}\, d\Sv \Big)\label{eq:var_prob_huwa_F_pairing}.
\end{align}
Above, $\uv_{\Tv}$ and $P_{\Nv\Nv}$ have to satisfy the corresponding boundary conditions on $\Gamma_D$ and $\Gamma_N$, 
respectively. The work of external forces is adapted according to \eqref{eq:Wtdnns}.

The three-field formulation \eqref{eq:var_prob_huwa_F} allows all kind of nonlinear material laws. As we use 
discontinuous elements for the additional deformation gradient field $\Ft$, it can be eliminated on element 
level, which makes the method competitive. Additionally, $\Pt$ can be eliminated by hybridization. Therefore, as in the linear case, a vector field $\alphav$ in
normal direction is added on all element interfaces. The first Piola--Kirchhoff stress is then assumed piecewise
smooth but discontinuous, and the hybridized optimization problem reads for all $\uv$ and $\alphav$ with
$\uv_{\Tv} =(\uv_D)_{\Tv}$ and $\alpha_{\Nv} = (u_D)_{\Nv}$ on $\Gamma_D$:
\begin{align}
\Lag^{\Ft}_h(\uv, \Ft,\Pt, \alphav) \to &\min_{\uv \text{ adm.}} \min_{\alphav \text{ adm.}} \max_{\Pt \text{ disc.}} \min_{\Ft \text{ disc.}}, \label{eq:laghF1}\\
\text{with } \Lag^{\Ft}_h(\uv, \Ft,\Pt, \alphav) &:= \Lag^{\Ft}(\uv, \Ft,\Pt) + \sum_{T \in \T} \int_{\partial T} P_{\Nv\Nv} \alpha_{\Nv}\, d\Sv,\label{eq:laghF2}
\end{align}
where in the definition of $\Lag^{\Ft}$ the work of external forces $W_{ext}^{\text{TDNNS}}$ is implicitly replaced by $W_{ext}^{\text{TDNNSh}}$ given in \eqref{eq:WTDNNSh}.

After static condensation of $\Ft$ and $\Pt$, a minimization problem in tangential-continuous $\uv$ and normal-continuous
$\alphav$ remains, see \cite{Neun21} for more details.

We note that the identity matrix $\It$ can exactly be represented by the discontinuous elements used for $\Ft$. Thus, it is equivalent to use $\Gt := \Ft -\It = \nabla u$ as independent field instead of $\Ft$.

\paragraph{Gradient splitting}
However, we can further simplify \eqref{eq:laghF1}--\eqref{eq:laghF2} by using the additive splitting of the deformation gradient into a symmetric and a
skew-symmetric part, $\Ft = \Ft_{\sym} + \Ft_{\skw}$. We note that the symmetric part $\Ft_{\sym}$ is related to the 
linearized strain tensor $\epst(\uv)$, while the skew-symmetric part $\Ft_{\skw}$ is equivalent to $\opcurl(\uv)$.
\begin{align}
\Ft(\uv) &= \Ft_{\sym}(\uv) + \Ft_{\skw}(\uv), & \Ft_{\sym}(\uv) &=  \It + \epst(\uv), & \Ft_{\skw}(\uv) &= \skw(\text{curl}(\uv)),\label{eq:add_spl_F}
\end{align}
with the $\skw$-operator defined as
\begin{align}
\skw(\vv) := \frac{1}{2}\begin{pmatrix}
\,0 & -\vv \\ \,\vv & \,0
\end{pmatrix} \text{ in 2D}, && \skw(\vv) := \frac{1}{2}\begin{pmatrix}
\,0 & -\vv_3 & \,\vv_2\\ 
\,\vv_3 & \,0 & -\vv_1\\
-\vv_2 & \,\vv_1 & \,0
\end{pmatrix} \text{ in 3D}.\label{eq:def_skw_op}
\end{align}
From theory for Maxwell's equations (see e.g. the monograph by Monk \cite{Mon03}), it is well-known that for
a piecewise smooth and tangentially continuous vector field, the $\opcurl$-operator is well-defined in
the sense of a square-integrable function. Thus, there is no need to ``lift'' the skew-symmetric part
$\Ft_{\skw}(\uv)$ to an independent $\Ft_{\skw}$, as $\Ft_{\skw}(\uv)$ is already in $L^2(\Omega)$.
The (hybridized) Lagrangian from \eqref{eq:laghF2} can be adapted accordingly, now using only the symmetric part
$\Pt_{\sym}$ of the first Piola--Kirchhoff stress as a multiplier for the constraint $\Ft_{\sym} = \Ft_{\sym}(\uv)$:
\begin{align}
\Lag^{\Ft,\sym}_h(\uv, \Ft_{\sym},\Pt_{\sym}, \alphav) \to &\min_{\uv \text{ adm.}} \min_{\alphav \text{ adm.}} \max_{\Pt_{\sym} \text{ disc.}} \min_{\Ft_{\sym} \text{ disc.}}, \label{eq:laghFsym1}\\
\text{with } \Lag^{\Ft,\sym}_h(\uv, \Ft_{\sym},\Pt_{\sym}, \alphav) :=  &\ 
\int_{\Omega}\Psi(\Ft_{\sym} + \skw(\text{curl}(\uv)))\,d\Xv - W_{ext}^{\text{TDNNSh}} \nonumber\\
& + \langle \epst (\uv), \Pt_{\sym}\rangle_\T - \int_\Omega (\Ft_{\sym}-\It):\Pt_{\sym}\, d\Xv\label{eq:laghFsym2}\\
& + \sum_{T \in \T} \int_{\partial T} P_{\sym,\Nv\Nv} \alpha_{\Nv}\, d\Sv. \nonumber
\end{align}

Note that the normal-normal continuity condition for the first Piola--Kirchhoff stress is equivalent
to the continuity of the normal-normal component of its symmetric part, as
\begin{align}
P_{\Nv\Nv} &= P_{\sym,\Nv\Nv} & \text{and}&& P_{\skw,\Nv\Nv} = 0.
\end{align}
In the definition of the Lagrangian \eqref{eq:laghFsym2}, we re-use the distributional definition of the duality pairing
$\langle \epst(\cdot), \cdot \rangle_\T$ from the linear TDNNS method in \eqref{eq:def_dual_pair2}.
Stress elements for the linear case can be used for the discretization of $\Pt_{\sym}$. 
Further details on the choices of finite elements, also for the symmetric deformation gradient
$\Ft_{\sym}$, are provided in Section~\ref{sec:numerics}. This simplification leads to fewer local
degrees of freedom than the original hybridized equation \eqref{eq:laghF1}--\eqref{eq:laghF2}.
One may recover $\Pt_{\skw}$ in a post-processing step by the equation $\frac{\partial\Psi}{\partial \skw(\text{curl}(\uv))}=\Pt_{\skw}$.

\subsection{Lifting to $\Ct$}
\label{subsec:lifting_C}
Due to objectivity of materials the energy potential $\Psi$ depends on the Cauchy--Green strain tensor $\Ct$. Thus, instead of solving the constrained minimization problem \eqref{eq:min_prob_huwa_F} we can postulate for a piecewise smooth and globally continuous $\tilde{\uv}$, satisfying $\tilde{\uv}=\uv_D$ on $\Gamma_D$, the problem
\begin{align}
\int_{\Omega}\Psi(\Ct)\,d\Xv -W_{ext}&\to \min\limits_{\substack{\tilde{\uv}\text{ cont.}\\\Ct=\Ct(\tilde{\uv})}},\label{eq:min_prob_huwa_C}
\end{align}
with the corresponding Lagrangian
\begin{align}
\mathcal{L}(\uv, \Ct,\Sigmat)&:=\int_{\Omega}\Psi(\Ct)\,d\Xv-W_{ext} - \int_{\Omega}\frac{1}{2}\Big(\Ct-\underbrace{(\nabla\uv+\It)^\top(\nabla\uv+\It)}_{=\Ct(\uv)}\Big):\Sigmat\,d\Xv.\label{eq:lag_func_huwa_cont_C}
\end{align}
Here, the Lagrange multiplier $\Sigmat$ is given by the second Piola--Kirchhoff stress tensor, which can be readily seen by taking the first variation of \eqref{eq:lag_func_huwa_cont_C} in direction $\delta\Ct$ together with \eqref{eq:pk2}
\begin{align}
\int_{\Omega}\frac{\partial\Psi}{\partial\Ct}(\Ct):\delta\Ct -\frac{1}{2}\Sigmat:\delta\Ct\,d\Xv = 0\qquad \forall \delta\Ct.
\end{align}
As discussed in Subsection \ref{subsec:lifting_F} the integral in \eqref{eq:lag_func_huwa_cont_C} is not well defined for tangential-continuous displacement fields. Further, the balance equation \eqref{eq:balance} implies that the first Piola--Kirchhoff tensor has to be normal-continuous rather than $\Sigmat$. Motivated by this, the second Piola--Kirchhoff stress tensor $\Sigmat$ and the Cauchy--Green strain tensor $\Ct$ are assumed to be piecewise smooth discontinuous and the hybridization variable $\alpha$ is used to enforce the normal-normal continuity of $\Pt=\Ft\Sigmat$ yielding the following saddle point problem
\begin{align}
\Lag_h^{\Ct}(\uv, \Ct,\Sigmat,\alphav)\to&\min_{\uv \text{ adm.}} \min_{\alphav \text{ adm.}} \max_{\Sigmat \text{ disc.}} \min_{\Ct \text{ disc.}},\label{eq:lag_C1} \\
\text{with }\Lag_h^{\Ct}(\uv, \Ct,\Sigmat,\alphav)&:=\int_{\Omega}\Psi(\Ct)\,d\Xv-W_{ext}^{\text{TDNNSh}} + \frac{1}{2}\langle\Ct(\uv),\Sigmat\rangle_{\T}- \frac{1}{2}\int_{\Omega}\Ct:\Sigmat\,d\Xv \nonumber\\
&\quad+ \sum_{T\in\T}\int_{\partial T}\left(\Ft(\uv)\Sigmat\right)_{\Nv\Nv}\alphav_{\Nv}\,d\Sv\label{eq:lag_C2}\\
\text{and } \langle \Ct(\uv), \Sigmat\rangle_{\T} &:=
\sum_{T \in \T} \Big( \int_T \Ct(\uv) : \Sigmat\,d\Xv - 2\int_{\partial T} \uv_{\Nv}( \Ft(\uv)\Sigmat)_{\Nv\Nv}\, d\Sv \Big).
\end{align}
The proof that this method is consistent, i.e., the smooth exact solution $\tilde{\uv}$ of \eqref{eq:balance} together with $\alphav:=\tilde{\uv}_{\Nv}\Nv$, $\Ct:=\Ct(\tilde{\uv})$, and $\Sigmat:=2\frac{\partial\Psi}{\partial\Ct}(\Ct(\tilde{\uv}))$  solves \eqref{eq:lag_C1} can be found in Appendix~\ref{app:consistency_C}.

Again, due to the discontinuity of $\Sigmat$ and $\Ct$, they can be eliminated at element level and the resulting system involving $\uv$ and $\alphav$ is positive. Therefore, the same amount of coupling dofs are needed as for the lifting of $\Ft$ in the previous subsection.

Linearizing \eqref{eq:lag_C2}, see Appendix \ref{app:linear_C}, yields that in the small deformation regime \eqref{eq:laghF1} and \eqref{eq:lag_C2} coincide. Further, assuming a quadratic potential $\Psi$ the hybridized TDNNS method \eqref{eq:var_prob_lin_tdnns_hyb} can be recovered.

We note that due to the affine relation $\Et = 0.5(\Ct-\It)$ and
\begin{align}
\frac{\partial \Psi}{\partial \Et}:\delta\Et = \frac{\partial\Psi}{\partial\Ct}:\delta\Ct\label{eq:rel_delta_E_C}
\end{align}
the problem
\begin{align}
\int_{\Omega}\Psi(\Et)\,d\Xv -W_{ext}\to \min\limits_{\substack{\tilde{\uv}\text{ cont.} \\\Et=\Et(\tilde{\uv})}}\label{eq:min_prob_huwa_E}
\end{align}
is equivalent to \eqref{eq:min_prob_huwa_C} and thus, one could use $\Et$ instead of $\Ct$ as additional field.

\paragraph{Stabilization techniques}
To improve robustness of this method in the large deformation regime, one may add the following well known stabilization term
\begin{align}
\sum_{T\in\T}\int_{\partial T}\frac{c_1}{h} (\uv-\alphav)_{\Nv}(\uv-\alphav)_{\Nv}\,d\Sv\label{eq:stab_u_alpha_n}
\end{align}
from Hybrid Discontinuous Galerkin (HDG) techniques \cite{BBF13,CGL09} to \eqref{eq:lag_C2}. Here, $h$ denotes the ratio of the element volume and the boundary area, $h=\frac{J}{J_{\text{bnd}}}$ ($J$ and $J_{\text{bnd}}$ are the volume and area measures of the mapping from the reference to the physical element) used especially for anisotropic elements and $c_1>0$ is a positive constant. As there holds for the true solution $\tilde{\uv}_{\Nv}=\alphav_{\Nv}$, \eqref{eq:stab_u_alpha_n} is consistent.

A second stabilization technique consists of adding the consistent term
\begin{align}
\sum_{T\in\T}\int_T c_2(\Ct-\Ct(u)):(\Ct-\Ct(u))\,d\Xv\label{eq:stab_C_Cu}
\end{align}
enforcing the element-wise equality of the lifting with $c_2>0$. Note, that these terms increase the stability but lead to less accurate solutions if the stability parameters $c_1$ and $c_2$ are chosen too large.

\subsection{Lifting to $\Ft$ and projection to $\Ct$}
\label{subsec:lifting_F_proj_C}
In numerical experiments we observed that the method presented in Subsection \ref{subsec:lifting_F}, where a lifting of $\Ft$ is considered, is more robust compared to the approach in the previous subsection, which may lead to more accurate solutions. The number of Newton iterations needed for one load step to converge, however, is significantly higher than for the second method. This motivates to combine both, leading to a lifting of $\Ft$, like in the first method. Then, the corresponding Cauchy--Green strain tensor $\Ct(\Ft) =\Ft^\top\Ft$ is going to be interpolated to an independent new field $\Ct$. This can be interpreted as a projection of $\Ct(\Ft)$ to $\Ct$ and thus, differs compared to the second approach, where $\Ct$ is a lifting rather than a projection.

Again, we start with a constraint minimization problem for a piecewise smooth, globally continuous and admissible $\tilde{\uv}$
\begin{align}
\int_{\Omega}\Psi(\Ct)\,d\Xv -W_{ext}\to \min\limits_{\substack{\tilde{\uv}\text{ cont.}\\\Ft=\Ft(\tilde{\uv})\\\Ct=\Ct(\Ft)}}
\end{align}
together with the Lagrangian
\begin{align}
\Lag(\uv,\Ft,\Pt,\Ct,\Sigmat):= \int_{\Omega}\Psi(\Ct)\,d\Xv -W_{ext} - \int_{\Omega}(\Ft-\nabla\uv-\It):\Pt\,d\Xv - \int_{\Omega}\frac{1}{2}(\Ct-\Ft^\top\Ft):\Sigmat\,d\Xv.\label{eq:lag_FC_cont}
\end{align}
Computing the variation in direction $\delta\Ct$ yields like in the previous subsection that the Lagrange multiplier $\Sigmat$ is the second Piola--Kirchhoff stress tensor. With the variation in direction $\delta \Ft$ we obtain that $\Pt$ corresponds to the first Piola--Kirchhoff stress tensor
\begin{align}
\int_{\Omega}\Pt:\delta\Ft-\Ft\Sigmat:\delta\Ft\,d\Xv = 0\qquad \forall \delta\Ft.
\end{align}
The integral involving the gradient of $\uv$ in \eqref{eq:lag_FC_cont} is interpreted as a distribution as in Subsection \ref{subsec:lifting_F} for piecewise smooth and tangential-continuous displacement fields. Together with $\Pt$ piecewise smooth, non-symmetric and $P_{\Nv\Nv}$ continuous, and $\Ft$, $\Sigmat$ and $\Ct$ piecewise smooth discontinuous we postulate the following saddle point problem
\begin{align}
\Lag^{\Ft\Ct}(\uv, \Ft,\Pt,\Ct,\Sigmat) \to &\min_{\uv \text{ adm.}} \max_{\Pt \text{ adm.}} \min_{\Ft \text{ disc.}}\min_{\Ct \text{ disc.}}\max_{\Sigmat \text{ disc.}},\label{eq:lag_FC1}\\
\text{with }\Lag^{\Ft\Ct}(\uv,\Ft,\Pt,\Ct,\Sigmat)&:= \int_{\Omega}\Psi(\Ct)\,d\Xv-W_{ext}^{\text{TDNNS}} - \langle \Ft-\nabla\uv-\It,\Pt\rangle_\T - \int_{\Omega}\frac{1}{2}(\Ct-\Ft^\top\Ft):\Sigmat\,d\Xv,\label{eq:lag_FC2}
\end{align}
where $\langle\cdot,\cdot\rangle_\T$ is defined as in \eqref{eq:var_prob_huwa_F_pairing}. 

One may pose the hybridized Lagrangian $\Lag^{\Ft\Ct}_h(\uv, \Ft,\Pt,\Ct,\Sigmat,\alphav)$ according to \eqref{eq:laghF1}--\eqref{eq:laghF2}. Further, in analogy to Subsection \ref{subsec:lifting_F} the gradient splitting can be used to eliminate the skew-symmetric part of $\Ft$ and $\Pt$ yielding the Lagrangian $\Lag^{\Ft\Ct,\sym}_h(\uv, \Ft_{\sym},\Pt_{\sym},\Ct,\Sigmat,\alphav)$.

The additional term
\begin{align}
\int_{\Omega}(\Ct-\Ft^\top\Ft):\Sigmat\,d\Xv\label{eq:dual_pairing_l2}
\end{align}
is well defined for piecewise smooth, discontinuous and thus, square-integrable functions. As a result $\Ct$ is the (local) $L^2$-projection of $\Ft^\top\Ft$ onto the space of polynomial order $k$ used for $\Sigmat$, see Section \ref{sec:numerics},
\begin{align}
\Ct = \mathcal{I}_{L^2}^k(\Ft^\top\Ft).
\end{align}
We note that this method is equivalent to the first approach if the polynomial order for $\Ct$ and $\Sigmat$ are chosen large enough, namely twice as the degree used for $\Ft$. Further, in the small strain regime $\Lag^{\Ft\Ct}_h(\uv, \Ft,\Pt,\Ct,\Sigmat,\alphav)$ reduces also to the hybridized linear TDNNS method \eqref{eq:var_prob_lin_tdnns_hyb}, as \eqref{eq:dual_pairing_l2} becomes
\begin{align}
\int_{\Omega}(\Ct-\It+\epst(\uv)):\Sigmat\,d\Xv.
\end{align}

For this method $\Sigmat$, $\Pt$, $\Ct$, and $\Ft$ can be eliminated by static condensation techniques at element level leading to a minimization problem in $\uv$ and $\alphav$. Thus, the number of coupling dofs coincide with those in Subsections \ref{subsec:lifting_F} and \ref{subsec:lifting_C}. However, the number of local dofs is higher compared to the previous methods.

\section{Updated Lagrangian}
\label{sec:updated_lagrangian}
In problems with large rotations, normal and tangential directions vary strongly when going from reference to actual configuration. This may lead to problems -- for the methods with lifting to $\Ft$ or $\Ft/\Ct$, we observed less than optimal behavior for thin structured elements, while for the method with lifting to $\Ct$ even breakdown was observed as the rotations approached $90^\circ$. To alleviate these problems, we propose to use an Updated Lagrangian scheme, where in each load step the configuration obtained in the last load step is used as intermediate configuration. As one may assume rotations within a single load step to be less than $90^\circ$, better behavior of the method is expected.

To describe the Updated Lagrangian based method, we follow \cite{P19} and introduce the notion of an \emph{intermediate configuration}. 
Starting from reference configuration $\Omega$, the deformation $\Xv \to \xv=\Xv + \uv(\Xv)$ associated to the displacement field $\uv$ maps 
$\Omega$ to the spatial configuration $\omega \subset \mathbb R^d$. In the following we introduce a multiplicative
splitting of the deformation, introducing some  ``intermediate'' displacement $\uv_\imt$ and a ``displacement update'' $\uv_\upd$,
\begin{align}
\uv(\Xv) &= \uv_\imt(\Xv) + \uv_\upd(\xv_0(\Xv))& \text{with} && \xv_\imt(\Xv) = \Xv + \uv_\imt(\Xv).
\end{align}
While the intermediate displacement $\uv_\imt$ is known and fixed in each load step, the displacement update $\uv_\upd$ is the unknown quantity to be computed.
\begin{figure}
	\centering
	\includegraphics[width=0.25\linewidth]{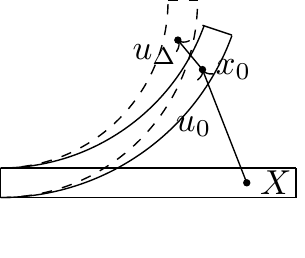}\hspace*{2cm}
	\includegraphics[width=0.25\linewidth]{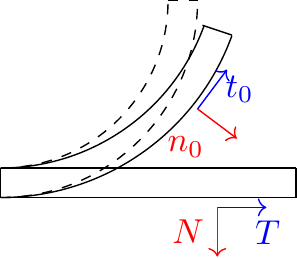}
	\caption{Reference configuration, intermediate configuration and actual configuration.}
	\label{fig:configs}
\end{figure}
The corresponding configuration $\omega_\imt := \xv_0(\Omega)$ shall be referred to as intermediate configuration. Boundary parts $\Gamma_D$ and $\Gamma_N$ are transformed to $\gamma_{D,\imt}$ and $\gamma_{N,\imt}$, respectively. For a graphical visualization of the different configurations see Figure~\ref{fig:configs}. One easily verifies the multiplicative splitting of the deformation gradient $\Ft$,
\begin{align}
\Ft(\uv) & = \Ft_\upd(\uv_\upd) \Ft_\imt(\uv_\imt)
\end{align}
with
\begin{align}
\Ft_\imt = \Ft_\imt(\uv_\imt) &= \It + \nabla_{\Xv}\uv_\imt, & \Ft_\upd(\uv_\upd) &= \It + \nabla_{\xv_\imt} \uv_\upd.
\end{align}
Above and in the following we will omit writing the explicit dependence of $\Ft_\imt$ on $\uv_\imt$, as $\Ft_\imt$ will not be used as an independent unknown. For the update gradient $\Ft_\upd(\uv_\upd)$ we explicitely state the dependence in all occurrances, to distinguish the gradient from the independent $\Ft_\upd$ introduced in the lifting procedure later.

The triangulation $\T$ of $\Omega$ is mapped to a -- possibly curved -- triangulation $\T_\imt$ of $\omega_\imt$. The unit outward normal vector $N$ in reference configuration is transformed to $\nv_\imt$ with respect to the intermediate configuration following
\begin{align}
\nv_\imt &= \frac{J_\imt}{J_{\text{bnd},\imt}} \Ft_{\imt}^{-\top} \Nv & \text{with}
&& J_\imt &= \det \Ft_\imt,\quad
J_{\text{bnd},\imt} = J_\imt  \|\Ft_{\imt}^{-\top} \Nv\|_2.
\end{align}
Note that in the definition above, $J_{A,\imt}$ denotes the transformation of the area element with normal~$\Nv$, while  $J_\imt$ is the transformation of the volume element. Normal and tangential components of vector or tensor fields are then defined using the transformed normal $\nv_\imt$ and denoted by subscripts $\nv_\imt$ and $\tv_\imt$, respectively.

Using the correct transformation rules, the deformation measures $\Ct$ and $J$ can be decomposed. To this end, let $\Ct_\upd(\uv_\upd) = (\Ft_\upd(\uv_\upd))^\top \Ft_\upd(\uv_\upd)$ be the Cauchy-Green tensor associated to $\uv_\upd$ acting on $\omega_\imt$. By $J_\upd(\uv_\upd) = \det(\Ft_\upd(\uv_\upd))$ we mean the according Jacobi determinant. They are connected to their absolute counterparts $\Ct(\uv)$ and $J(\uv)$ via
\begin{align}
\Ct(\uv) &= \Ft_\imt^\top \Ct_\upd(\uv_\upd) \Ft_\imt, & J(\uv)  = J_\imt J_\upd(\uv_\upd).
\end{align}

As mentioned before, the intermediate configuration $\omega_0$ and the corresponding displacement part $\uv_\imt$ is assumed known and constant in each load step. The displacement update $\uv_\upd$ on the other hand is a-priori unknown and a finite element solution shall be found. In this setting, the additional kinematic unknowns to be introduced will be independent $\Ft_\upd$ and $\Ct_\upd$ rather than $\Ft$ or $\Ct$.

The internal energy, which appears in all formulations, can then be re-written after a transformation of the volume integral using the independent $\Ft_\upd$ or $\Ct_\upd$
\begin{align}
\int_\Omega \Psi(\Ft)\, d\Xv &= \int_{\omega_\imt} J_\imt^{-1} \Psi( \Ft_\upd \Ft_\imt)\, d\xv_\imt,\\
\int_\Omega \Psi(\Ct)\, d\Xv &= \int_{\omega_\imt} J_\imt^{-1} \Psi(\Ft_\imt^\top \Ct_\upd \Ft_\imt)\, d\xv_\imt.
\end{align}

\subsection{Lifting to $\Ft_\upd$}
For a moment, assume the displacement function weakly differentiable, $\uv_\upd \in [H^1(\omega_\imt)]^d$. The correct Lagrange functional corresponding to \eqref{eq:lag_huwa_F} using the intermediate configuration and multiplicative decomposition of displacements reads
\begin{align}
\mathcal L(\uv_\upd, \Ft_\upd, \Pt_\imt) &:= \int_{\omega_\imt} J_\imt^{-1} \Psi( \Ft_\upd \Ft_\imt)\, d\xv_\imt - W_{ext,\imt} - \int_{\omega_\imt} ( \Ft_\upd - \It - \nabla_{\xv_\imt} \uv_\upd):\Pt_\imt\, d\xv_\imt,
\end{align}
with the work of external forces
\begin{align}
W_{ext,\imt} &= \int_{\omega_\imt} J_\imt^{-1} \fv \cdot \uv_\upd\, d\xv_\imt + \int_{\gamma_{N,\imt}} J_{\text{bnd},\imt}^{-1} \,\gv \cdot \uv_\upd\, ds_\imt.
\end{align}
One easily verifies that the Lagrangian multiplier $\Pt_\imt$ is connected to the Piola-Kirchhoff stress $\Pt$ via
\begin{align}
\Pt_\imt &= J_\imt^{-1} \Pt \Ft_{\imt}^\top.
\end{align}

As in Subsection~\ref{subsec:lifting_F}, we can now derive a (hybridized) TDNNS method in intermediate configuration. For $\uv_\upd$ piecewise smooth and tangential-continuous on the triangulation $\T_\imt$ of $\omega_\imt$, $\Pt_\imt$ piecewise smooth non-symmetric with $(P_\imt)_{\nv_\imt\nv_\imt}$ continuous, and $\Ft_\upd$ piecewise discontinuous we obtain a saddle point problem similar to \eqref{eq:var_prob_huwa_F}--\eqref{eq:var_prob_huwa_F_pairing},
\begin{align}
\mathcal L^{\Ft}_\imt(\uv_\upd,\Ft_\upd,\Pt_\imt) \to &\min_{\uv_\upd \text{ adm.}} \max_{\Pt_\imt \text{ adm.}} \min_{\Ft_\upd \text{ adm.}},\label{eq:updlaghF1}\\
\begin{split}
	\text{with } \mathcal L^{\Ft}_\imt(\uv_\upd,\Ft_\upd,\Pt_\imt) \ := & \int_{\omega_\imt} J_\imt^{-1} \Psi( \Ft_\upd \Ft_\imt)\, d\xv_\imt - W_{ext,\imt}\\
	& + \langle \nabla_{\xv_\imt} \uv_\upd, \Pt_\imt\rangle_{\T_\imt} - \int_{\omega_\imt} (\Ft_\upd - \It):\Pt_\imt\, d\xv_0,
\end{split}\\
\begin{split}
W^{\text{TDNNS}}_{ext,\imt} =& \int_{\omega_\imt} J_\imt^{-1} \fv\cdot\uv_\upd\, d\xv_\imt
 + \int_{\gamma_{N,\imt}} J_{\text{bnd},\imt}^{-1}\, \gv_{\tv_\imt}\cdot (\uv_\upd)_{\tv_\imt}\, ds_\imt \\
 &- \int_{\gamma_{D,\imt}} (\uv_D-\uv_\imt)_{\nv_0} (P_\imt)_{\nv_0\nv_0}\, d\sv_\imt,
 \end{split}\\
\text{and } \langle \nabla_{\xv_\imt} \uv_\upd, \Pt_\imt\rangle_{\T_\imt} =& \sum_{T_\imt \in \T_\imt} \left(
\int_{T_\imt} \nabla_{\xv_\imt} \uv_\upd : \Pt_\imt\, d\xv_\imt 
- \int_{\partial T_\imt} (\uv_\upd)_{\nv_\imt}\, (P_\imt)_{\nv_\imt\nv_\imt}\, ds_{\imt} \right). \label{eq:updlaghF4}
\end{align}
Additionally, in \eqref{eq:updlaghF1}--\eqref{eq:updlaghF4} the following boundary conditions have to be satisfied for $\uv_\upd$ and $\Pt_\imt$ to  be admissible,
\begin{align}
(\uv_\upd)_{\tv_\imt} &= (\uv_D - \uv_\imt)_{\tv_\imt} & \text{and} &&
(P_\imt)_{\nv_\imt\nv_\imt} &= J_{A,\imt}^{-1} g_{\nv_0}.
\label{eq:updlag_bc}
\end{align}
As in \eqref{eq:laghF1}--\eqref{eq:laghF2} hybridization is again possible in order to regain a symmetric positive system. The hybridized problem reads for all $\uv_\upd$ and $\Ft_\upd$ as above, $\Pt_\imt$ piecewise smooth discontinuous not necessarily satisfying the boundary condition from \eqref{eq:updlag_bc}, and $\nv_0$-continuous Lagrangian multipliers $\alpha$ on element interfaces with $\alpha_{\nv_0} = (\uv_D - \uv_\imt)_{\nv_\imt}$ on $\gamma_{D,\imt}$:
\begin{align}
 (\uv_\upd, \Ft_\upd, \Pt_\imt, \alpha) &\to \min_{\uv_\upd \text{ adm.}}
\min_{\alpha \text{ adm.}} \max_{\Pt_\imt \text{ disc.}} \min_{\Ft_\upd \text{ adm.}},\\
\text{with } \Lag^{\Ft}_{\imt,h} &:= \Lag^{\Ft}_{\imt}(\uv_\upd, \Ft_\upd, \Pt_\imt)
+ \sum_{T_\imt \in \T_\imt} \int_{\partial T_\imt} (P_\imt)_{\nv_\imt\nv_\imt} \alpha_{\nv_\imt}\, d\sv_\imt,\\
&\text{ in which } W_{ext,\imt}^{\text{TDNNS}} \text{ is replaced by }\nonumber\\
W_{ext,\imt}^{\text{TDNNSh}}&
:= \int_{\omega_\imt} J_\imt^{-1} \fv\cdot\uv_\upd\, d\xv_\imt
+ \int_{\gamma_{N,\imt}} J_{\text{bnd},\imt}^{-1}\, (\gv_{\tv_\imt}\cdot (\uv_\upd)_{\tv_\imt}
+ g_{\nv_\imt} \, \alpha_{\nv_\imt})\, ds_\imt .
\end{align}

\subsection{Lifting to $\Ct_\upd$ or lifting to $\Ft_\upd$ with projection to $\Ct_\upd$}

The variants described in Subsection~\ref{subsec:lifting_C} as well as Subsection~\ref{subsec:lifting_F_proj_C} are transformed to the Updated Lagrangian setting in the same way. The update tensors $\Ct_\upd$ and, in latter case, also $\Ft_\upd$ are discretized independently.

\section{Numerics}
\label{sec:numerics}

\subsection{Finite element spaces}
\label{subsec:fe_spaces}

For the sake of simplicity of presentation, we assume $\T$ to be a simplicial triangulation. Note, however, that hybrid meshes can directly be treated using different elements in the same mesh. On the triangulation $\T$ we define the set of all piecewise polynomials up to order $k$ by $\Pi^k(\T)$ and the polynomials living only on the skeleton $\E$ are denoted by $\Pi^k(\E)$.

We define the following finite element spaces of order $k$
\begin{subequations}
\begin{align}
& V^k := \{\vv\in[\Pi^k(\T)]^d: \vv\text{ continuous}\},\label{eq:fe_space_h1}\\
& U^k := \{\uv\in[\Pi^k(\T)]^d:  \uv_{\Tv}\text{ continuous}\},\label{eq:fe_space_hcurl}\\
& \Gamma^k := \{\alphav\in[\Pi^k(\E)]^d: \alphav_{\Nv}\text{ continuous}\},\label{eq:fe_space_hybrid}\\
& \Sigma^k := \{\sigmat\in[\Pi^k(\T)]^{d\times d}_{\sym}: \sigma_{\Nv\Nv}\text{ continuous}\},\\
& \mathcal{R}^k := \{\Ct\in[\Pi^k(\T)]^{d\times d}_{\sym}: \Ct_{\Tv\Tv}\text{ continuous}\}.\label{eq:fe_space_hcurlcurl}
\end{align}

\begin{figure}[h!]
	\centering
	\begin{tabular}{cccc}
		\includegraphics[width=0.2\linewidth]{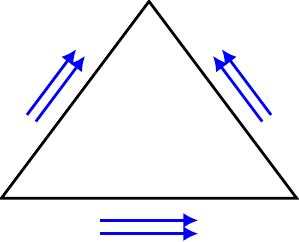} &
		\includegraphics[width=0.2\linewidth]{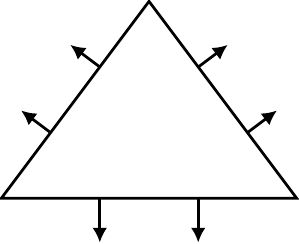} &
		\includegraphics[width=0.2\linewidth]{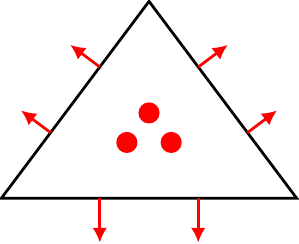} &
		\includegraphics[width=0.2\linewidth]{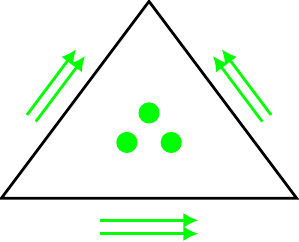}
	\end{tabular}
	\caption{N{\'e}d{\'e}lec, hybridization, stress and Regge finite element for $k=1$.}
	\label{fig:fe_spaces}
\end{figure}

\end{subequations}
The set of all piecewise polynomials which are globally continuous is denoted by $V$ and corresponds to the well known Lagrangian nodal finite elements \cite{ZT20}. The space $U$ is given by the N{\'e}d{\'e}lec elements, where the vector valued polynomials are tangential continuous over elements \cite{Nedelec1986}. The hybridization space $\Gamma$ can be implemented by using a facet space equipped with the normal vector \cite{BBF13}. As it follows the same transformation rules as H(div)-conforming Raviart--Thomas (RT) \cite{RT77} or Brezzi--Douglas--Marini (BDM) \cite{BDM85} elements, they can also be used (neglecting their inner dofs). The stress space $\Sigma$ consists of the normal-normal continuous elements introduced in \cite{PS11}. For two-dimensional domains, triangular and quadrilateral elements have been introduced, while for three-dimensional meshes tetrahedra, hexahedra and prismatric elements have been developed so far, see \cite{PS11,PS12,MPH:18}. The so-called Regge finite element space is constructed such that the tangential-tangential part $\Ct_{\Tv\Tv}:=(\It-\Nv\boldsymbol{\otimes}\Nv)\Ct(\It-\Nv\boldsymbol{\otimes}\Nv)$, one component in 2D and four in 3D, is continuous. Here $\boldsymbol{\otimes}$ denotes the dyadic product of two vectors. For a construction of Regge elements for triangles and tetrahedral we refer to \cite{Li18}. The construction of arbitrary order quadrilateral, prismatic, and hexahedra Regge elements is topic of an upcoming paper, see also \cite{Neun21}. Constructions based on hierarchical elements for the spaces \eqref{eq:fe_space_h1}--\eqref{eq:fe_space_hybrid} are given e.g. in \cite{Zaglmayr06}. In Figure \ref{fig:fe_spaces} triangular finite elements of \eqref{eq:fe_space_hcurl}--\eqref{eq:fe_space_hcurlcurl} of polynomial order $k=1$ are depicted.

One may break the continuity of the spaces $\Sigma^{k}$ and $\mathcal{R}^{k}$ for the presented methods, which will be denoted by $\Sigma^{k,\star}$ and $\mathcal{R}^{k,\star}$, respectively. Therefore, the same basis as for the continuous spaces are used, but no continuity over interfaces is required, i.e., the corresponding functions are discontinuous.\\

For the lifting of $\Ft$ method in Subsection \ref{subsec:lifting_F} the tangential-continuous displacement field $\uv$ will be discretized with N{\'e}d{\'e}lec elements $\uv\in U^k$, whereas for the normal-continuous hybridization variable $\alphav$ the space $\Gamma^k$ is used. For the Lagrangian $\Lag^{\Ft,\sym}_h(\uv, \Ft_{\sym},\Pt_{\sym}, \alphav)$ given in \eqref{eq:laghFsym2} the (symmetric part of the) first Piola--Kirchhoff stress tensor $\Pt_{\sym}$ is discretized by the stress space $\Sigma^{k,\star}$ such that the duality product $\langle\epst(\uv),\Pt_{\sym}\rangle_{\T}$ is well defined \cite{PS11}. Let us shortly motivate the choice of Regge-elements for the lifted deformation gradient $\Ft$: for the deformation gradient $\Ft$ of a piecewise smooth and globally continuous displacement field $\tilde{u}$ we observe that its gradient is tangential continuous, i.e., $\Ft_{\Tv}:=\Ft(\It-\Nv\boldsymbol{\otimes}\Nv)$ is continuous. Taking the symmetric part of $\Ft$ gives a tangential-tangential continuous field $\Ft_{\sym}$, motivating to use the discontinuous Regge finite elements, $\Ft_{\sym}\in\mathcal{R}^{k,\star}$. Further, the N{\'e}d{\'e}lec elements transform covariantly. More precisely, for the reference element $\hat{T}$ assume an affine transformation $\phi$ to a physical element, $\phi:\hat{T}\to T$, and define the corresponding gradient $\Gt=\nabla_{\hat{x}}\phi$. Let $\hat{\uv}$ be a N{\'e}d{\'e}lec element on the reference element $\hat{T}$. Then, with the definition $\uv\circ\phi = \Gt^{-\top}\hat{\uv}$, $\uv$ is a N{\'e}d{\'e}lec element on the physical element $T$, cf. \cite{Mon03}. The gradient of this tangential-continuous displacement field $\uv$ transforms with $(\nabla_\Xv \uv)\circ\phi=\Gt^{-\top}\nabla_{\hat{x}}\hat{\uv}\Gt^{-1}$, i.e., doubled covariantly. The tangential-tangential continuous Regge elements are constructed such that the mapping is also a double covariant transformation, $\Ft\circ\phi=\Gt^{-\top}\hat{\Ft}\Gt^{-1}$. Therefore, in $\int_{\Omega}(\Ft-\nabla \uv):\Pt\,d\Xv$ both, $\Ft$ and $\nabla \uv$, are transformed in the same manner from reference to affine-linear element. For non-affine transformations there holds (gradient of $\Gt$ has to be interpreted in the correct way) $(\nabla_{\Xv}u)\circ\Phi=\nabla_{\hat{x}}(\Gt^{-\top}\hat{u})\Gt^{-1}= \Gt^{-\top}(\nabla_{\hat{x}}\hat{u}-\nabla_{\hat{x}}\Gt^\top\Gt^{-\top}\hat{u})\Gt^{-1}$ and thus transforms also doubled covariantly. Summing up, for the first method \eqref{eq:laghFsym2} we seek for $(\uv,\alphav,\Ft_{\sym},\Pt_{\sym})$ in $U^k\times\Gamma^k\times \mathcal{R}^{k,\star}\times\Sigma^{k,\star}$.

The second approach uses the same spaces as the first one, we seek for $(\uv,\alphav,\Ct,\Sigmat)$ in $U^k\times\Gamma^k\times \mathcal{R}^{k,\star}\times\Sigma^{k,\star}$. Again, for a continuous $\tilde{u}$ it follows that $\Ft$ is tangential-continuous and thus, $\Ct=\Ft^\top\Ft$ is tangential-tangential continuous.

For the third method the solution $(\uv,\alphav,\Ft_{\sym},\Pt_{\sym},\Ct,\Sigmat)$ lives in $U^k\times\Gamma^k\times \mathcal{R}^{k,\star}\times\Sigma^{k,\star}\times \mathcal{R}^{k,\star}\times \mathcal{R}^{k,\star}$. The choices for $,\Ft_{\sym},\Pt_{\sym},\Ct$ follow the same ideas as before. Motivated by the $L^2$-interpolation term \eqref{eq:dual_pairing_l2} the same space for $\Sigmat$ is used as for $\Ct$.

\subsection{Examples}
All methods are implemented in the open source finite element library Netgen/NGSolve\footnote{www.ngsolve.org} \cite{Sch97,Sch14}. In what follows we will denote the first, second, and third method as  $\Ft$-, $\Ct$-, and $\Ft\Ct$-based method, respectively. We observed that the $\Ft$-based method is more robust compared to the $\Ct$-based method giving slightly better results, however more Newton iteration steps are needed for every load step. The $\Ft\Ct$-based method combines both advantages, the accuracy of $\Ft$- and the faster Newton convergence of $\Ct$-based method at the expense of being locally more expensive, due to the higher amount of local dofs.

For all benchmarks, if nothing stated otherwise, second order methods, $k=2$ are used. We stress that also polynomial degree $k=1$ can be used as well as higher polynomial orders. The methods are also compared with standard Lagrangian elements of polynomial order $k=2$ for the displacement $\uv$ denoted by method ``std''. We note that for the same grids the coupling dofs of our methods are nearly doubled compare to the standard Lagrangian method in two dimensions as there are asymptotically more edges than vertices, namely $\#E\approx 3\#V$, where the dofs are placed (see Figure~\ref{fig:fe_spaces}). In the three dimensional case the coupling type dofs are approximately four times due to the fact that $\#E\approx 7\#V$.

To solve the arising nonlinear problems a (damped) Newton method is used together with load steps to increase the right-hand side. As we started with (constraint) minimization problems the stiffness matrices appearing in the Newton iterations are symmetric. Further, due to the static condensation, the resulting (smaller) system is also a minimization problem, enabling the use of the built in \emph{sparsecholesky} solver. Full code examples are available\footnote{www.gitlab.com/mneunteufel/nonlinear\_elasticity}.

We will consider two different material laws of neo-Hooke as hyperelastic potentials, namely
\begin{subequations}
\begin{align}
\begin{cases}
\Psi_1(\Ft):=\frac{1}{2}\left(\text{tr}(\Ft^\top\Ft)-d\right) -\mu\log\det\Ft + \frac{\lambda}{2}(\det\Ft-1)^2,\\
\Psi_1(\Ct):=\frac{1}{2}\left(\text{tr}(\Ct)-d\right) -\frac{\mu}{2}\log\det\Ct + \frac{\lambda}{2}(\sqrt{\det\Ct}-1)^2,\\
\end{cases}\label{eq:neohooke_nolog}\\
\begin{cases}
\Psi_2(\Ft):=\frac{1}{2}\left(\text{tr}(\Ft^\top\Ft)-d\right) -\mu\log\det\Ft + \frac{\lambda}{2}(\log\det\Ft)^2,\\
\Psi_2(\Ct):=\frac{1}{2}\left(\text{tr}(\Ct)-d\right) -\frac{\mu}{2}\log\det\Ct + \frac{\lambda}{8}(\log\det\Ct)^2,
\end{cases}\label{eq:neohooke_log}
\end{align}
\end{subequations}

where $\mu$ and $\lambda$ are the Lam\'{e} parameters and $d=2,3$ the spatial dimension.

\subsubsection{Shearing Plate}
\label{subsubsec:shearing_plate}

\begin{figure}[h!]
	\centering
	\includegraphics[width=0.1\linewidth]{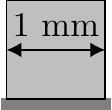}
	\caption{Geometry of shearing plate example.}
	\label{fig:shearing_plate_geo}
\end{figure}

\begin{figure}[h!]
	\centering
	\includegraphics[width=0.48\linewidth]{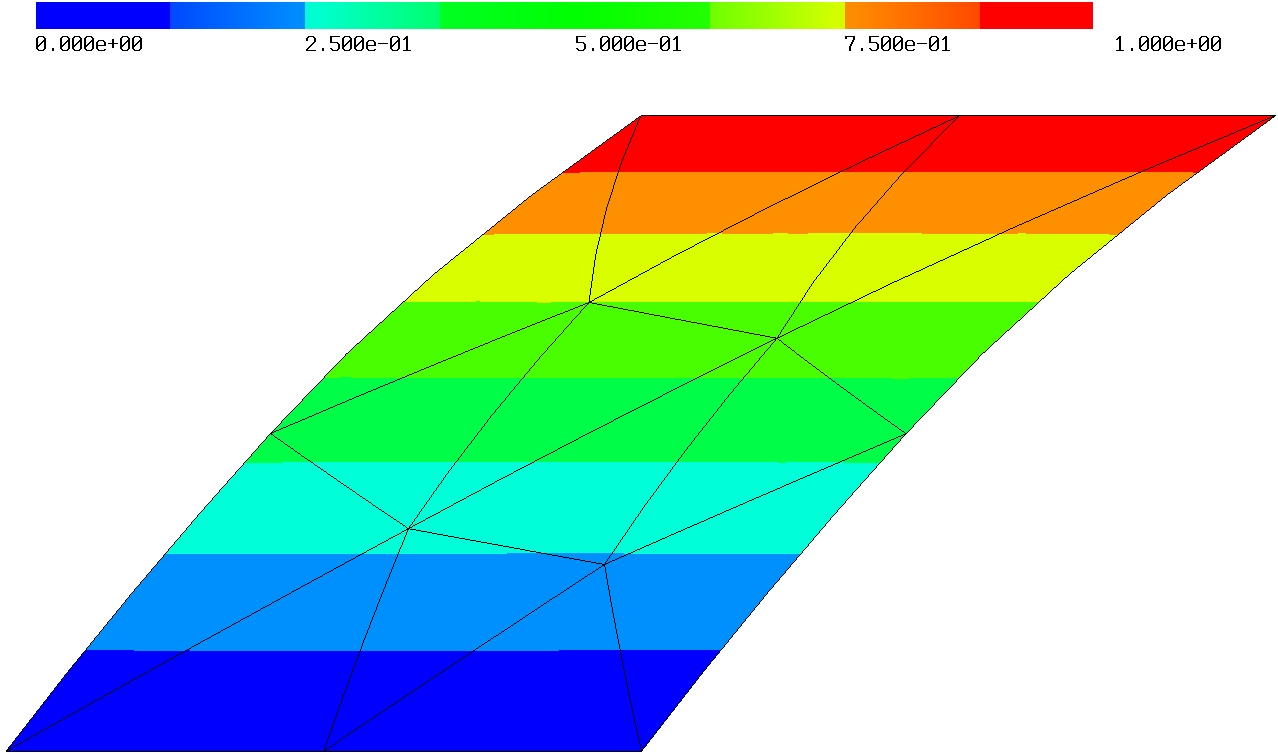}
	\includegraphics[width=0.48\linewidth]{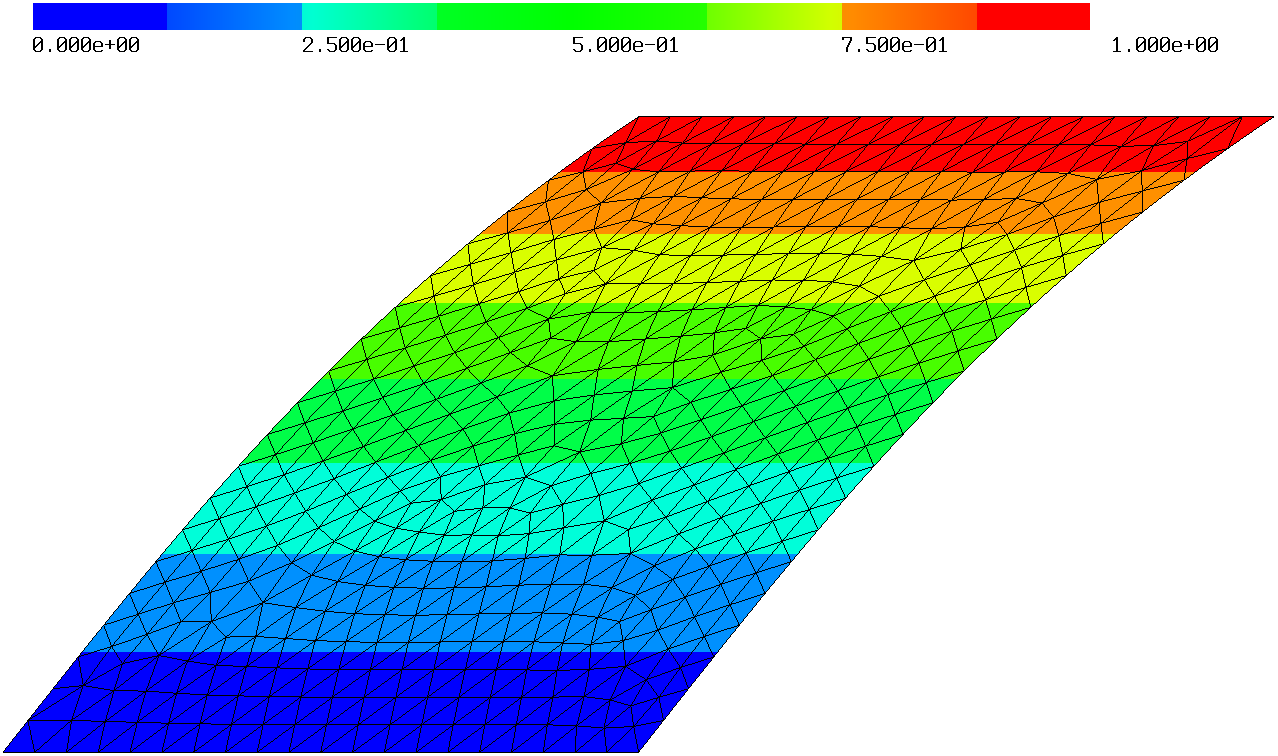}
	\caption{Final deformation of shearing plate example with 14 and 934 elements.}
	\label{fig:shearing_plate_displ}
\end{figure}

\begin{figure}[h!]
	\centering
	\includegraphics[width=0.45\linewidth]{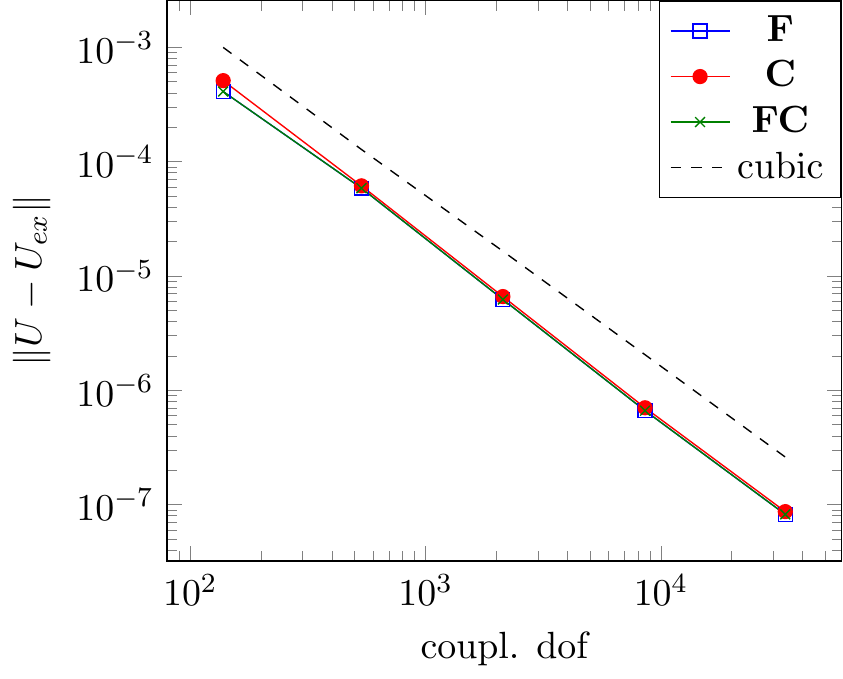}
	\includegraphics[width=0.45\linewidth]{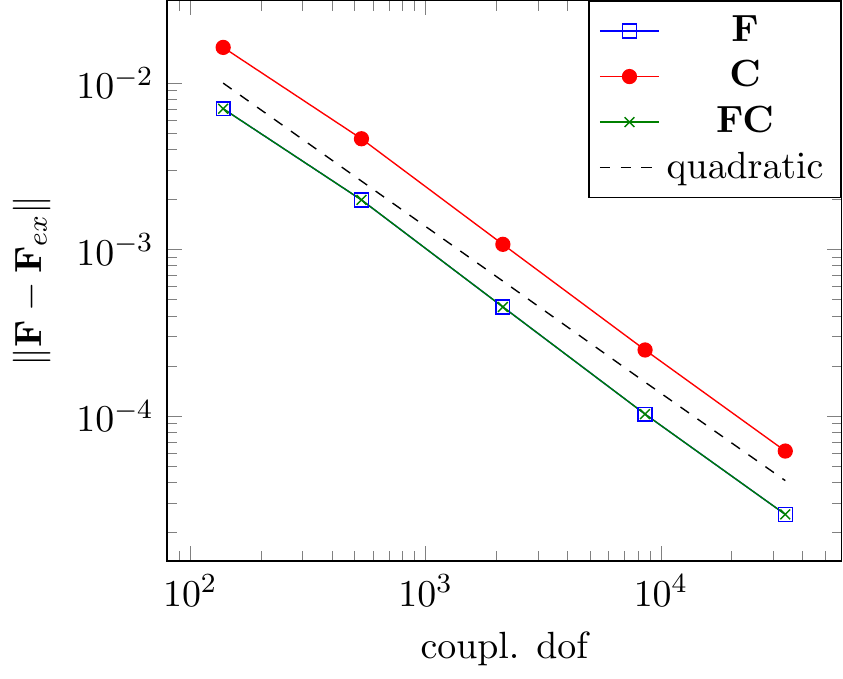}
	\caption{$L^2$-errors of the methods with respect to the coupling dofs for shearing plate example. Left: displacement error $\|U-U_{ex}\|$. Right: error of deformation gradient $\|\Ft-\Ft_{ex}\|$ for $\Ft$- and $\Ft\Ct$-based method and Cauchy--Green strain tensor $\|\Ct-\Ct_{ex}\|$ for $\Ct$-based method.}
	\label{fig:shearing_plate_res}
\end{figure}

A clamped square plate with length $\SI{1}{\milli\meter}$ is subjected to shear loads, see Figure \ref{fig:shearing_plate_geo}. This benchmark has been considered in \cite{ASY17,R02}. We assume the displacement field
\begin{align*}
U_{ex}=\begin{pmatrix}
\frac{1}{2}y^3 +\frac{1}{2}\sin\left(\frac{\pi\,y}{2}\right)\\ 0
\end{pmatrix}
\end{align*}
and the neo-Hookean material law \eqref{eq:neohooke_nolog} with parameters $\mu=\lambda=\SI{1}{\newton\per\milli\meter\squared}$ is used. With $U_{ex}$ at hand the corresponding right hand sides $f$ and $g$ can be readily computed. Unstructured triangular meshes are used and in Figure \ref{fig:shearing_plate_displ} the final deformation is depicted. The results for all methods can be found in Table \ref{tab:shearing_plate_res}, and Figure \ref{fig:shearing_plate_res} shows the error of the displacement $\|U-U_{ex}\|$ and deformation gradient $\|\Ft-\Ft_{ex}\|$. Note that for the $\Ct$-based method $\|\Ct-\Ct_{ex}\|$ is computed instead, as the approach does not use the deformation gradient $\Ft$. Thus, the corresponding error curve is shifted above compared to the others. We observe the expected optimal rates for all methods, cubic for the displacement $U$ and quadratic for $\Ft$ and $\Ct$, respectively. For the $\Ct$-based method, however, we needed stabilization \eqref{eq:stab_u_alpha_n} with $c_1=1$ to guarantee convergence for finer grids.

\begin{table}[h!]
\centering
\input{shearing_plate.out}
\caption{Results for shearing plate example. For all methods the number of elements, number of dofs and coupling dofs, and the $L^2$-error of the displacement and deformation gradient are presented. $^\star$: For the $\Ct$-based method $\|\Ct-\Ct_{ex}\|_{L^2}$ is computed.}
\label{tab:shearing_plate_res}
\end{table}

\subsubsection{Cook's Membrane}
\label{subsubsec:cooks_membrane}
\begin{figure}[h!]
	\centering
	\includegraphics[width=0.45\linewidth]{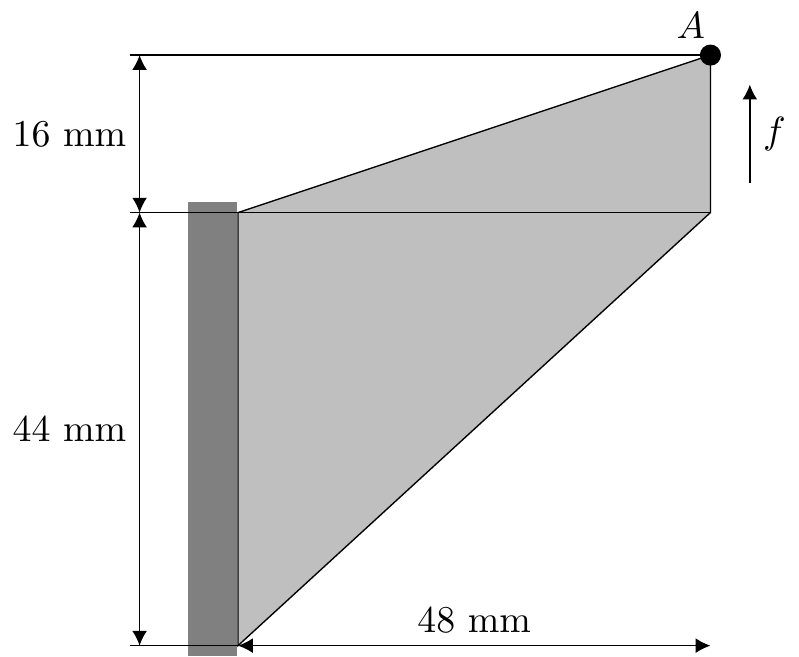}
	\caption{Geometry of Cook's membrane example.}
	\label{fig:cooks_membrane_geo}
\end{figure}

\begin{figure}[h!]
	\centering
	\includegraphics[width=0.6\linewidth]{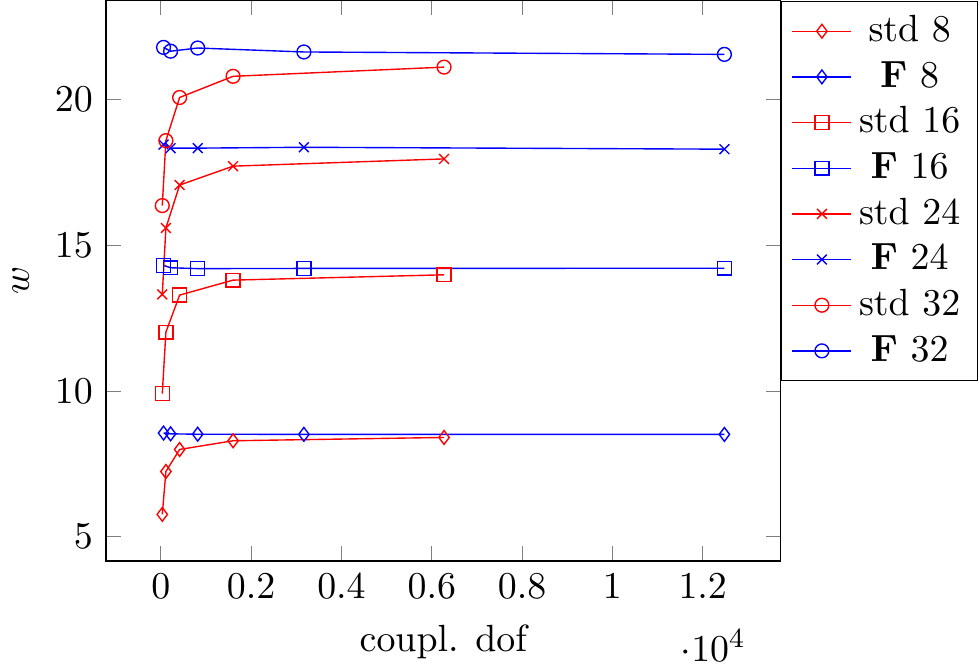}
	\caption{Vertical deflection $w$ at point $A$ for Cook's membrane example for $f=8,16,24,32$ with standard and $\Ft$-based method.}
	\label{fig:result_cooks_membrane}
\end{figure}

\begin{figure}[h!]
	\centering
	\includegraphics[width=0.48\linewidth]{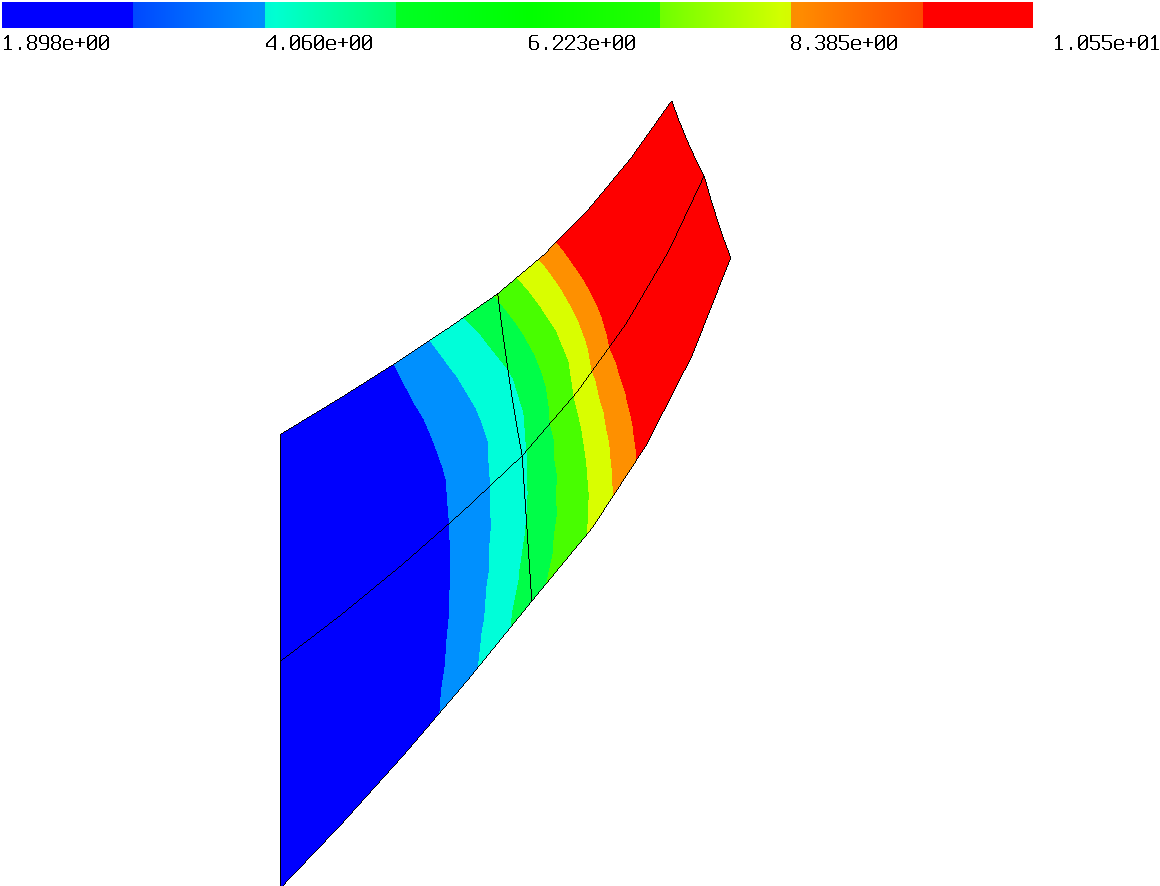}
	\includegraphics[width=0.48\linewidth]{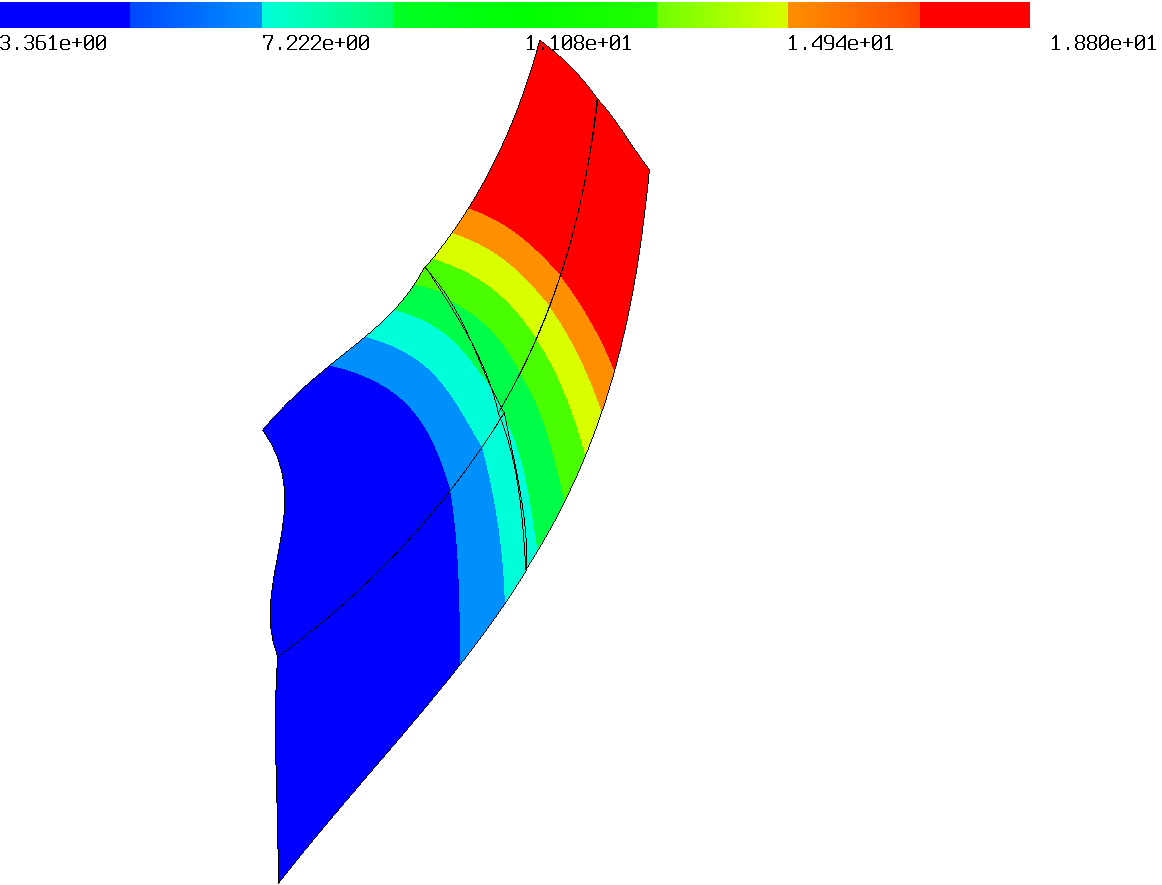}
	
	\includegraphics[width=0.48\linewidth]{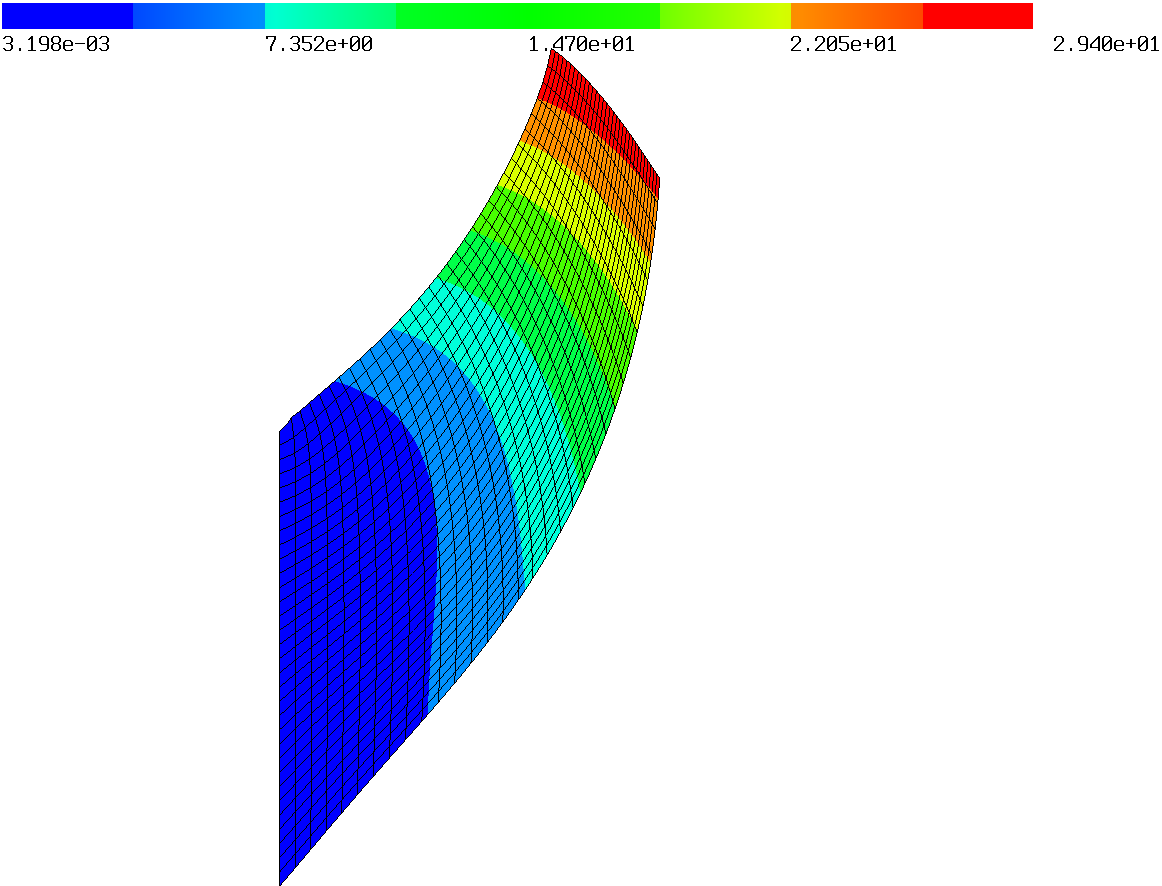}
	\includegraphics[width=0.48\linewidth]{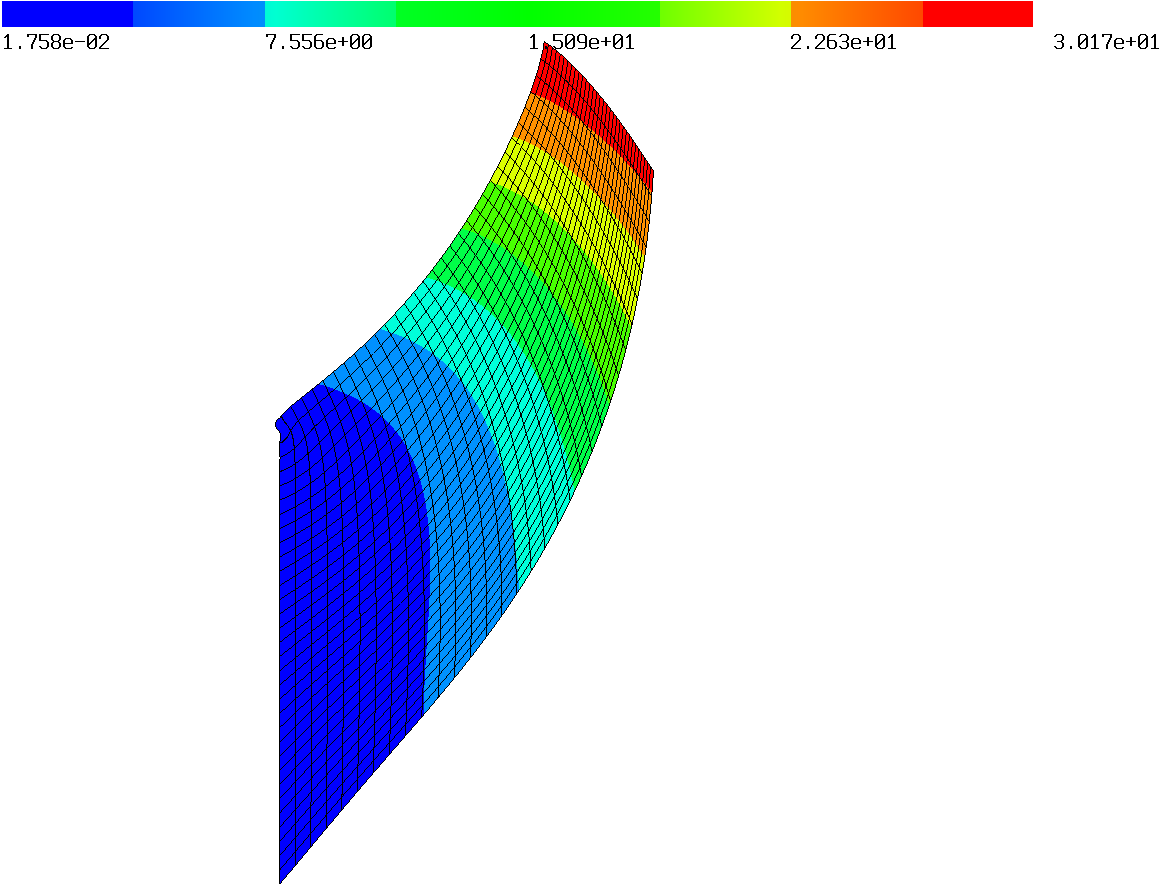}
	\caption{Final deformation for Cook's membrane example with $f=32$ standard method (left) and $\Ft$-based method (right) for $2\times2$ (top) and $32\times 32$ (bottom) grid.}
	\label{fig:cooks_membrane_displ_f32}
\end{figure}

We consider the Cook's membrane problem, see Figure \ref{fig:cooks_membrane_geo}, which has been used as a benchmark problem by \cite{ASY17,R02}. Material parameters for the hyperelastic potential \eqref{eq:neohooke_log} are $\mu= \SI{80.194}{\newton\per\milli\meter\squared}$ and $\lambda=\SI{40889.8}{\newton\per\milli\meter\squared}$, which results in nearly incompressible behavior. The quantity of interest is given by the vertical deflection at point $A$, cf. Figure \ref{fig:cooks_membrane_geo}. Different shear forces $f=8,16,24,\SI{32}{\newton\per\milli\meter\squared}$ are considered and structured quadrilateral meshes with $2\times2$, $4\times4$, $8\times8$, $16\times16$, and $32\times32$ grids are used. It is well known that on the top left corner a singularity leads to reduced convergence rates. For this reason mostly adaptive (triangular) meshes are used to resolve the singularity. For the proposed methods, however, already a coarse $2\times2$ quadrilateral grid produces surprisingly accurate results being already in the correct magnitude, see Figure \ref{fig:result_cooks_membrane}. In Figure~\ref{fig:cooks_membrane_displ_f32} a comparison between the standard and $\Ft$-based method for $f=32$ is shown, where a clear difference in the vertical deflection can be seen. We observed that the quadrilateral on the top left deforms also on the clamped boundary, as the components are not prescribed point-wise giving the proposed methods more flexibility.

For the $\Ct$-based method stabilization techniques \eqref{eq:stab_u_alpha_n} and \eqref{eq:stab_C_Cu} are used with parameters $c_1=c_2=\frac{\mu}{2}$. With them, only on the finest grid for the large forces $f=24,32$ Newton's method did not converge. The results for all forces can be found in Tables \ref{tab:cooks_membrane_result_f816} -- \ref{tab:cooks_membrane_result_f2432}. All results agree with those in \cite{ASY17,R02}.

\begin{table}[h!]
	\centering
	\begin{tabular}{cc|cc|cc}
		& coupl. dof & w & $\|U\|$& w & $\|U\|$ \\
		\hline
		&&$f=8$&&$f=16$&\\
		std & 32 & 5.767 & 86.168 & 9.902 & 153.983 \\
		& 112 & 7.236 & 113.156 & 12.012 & 199.393 \\
		& 416 & 7.992 & 128.938 & 13.281 & 232.024 \\
		& 1600 & 8.290 & 135.377 & 13.794 & 246.522 \\
		& 6272 & 8.403 & 137.821 & 13.977 & 252.014 \\
		\hline
		$\Ft$ & 60 & 8.554 & 141.578 & 14.299 & 261.363 \\
		&  216 & 8.527 & 140.770 & 14.218 & 259.260 \\
		& 816 & 8.514 & 140.441 & 14.188 & 258.691 \\
		& 3168 & 8.509 & 140.312 & 14.195 & 259.070 \\
		& 12480 & 8.507 & 140.282 & 14.201 & 259.345 \\
		\hline
		$\Ct$  & 60 & 8.219 & 134.345 & 13.854 & 248.796 \\
		& 216 & 8.377 & 137.746 & 14.026 & 253.700 \\
		& 816 & 8.442 & 139.008 & 14.092 & 255.674 \\
		& 3168 & 8.474 & 139.567 & 14.122 & 256.684 \\
		& 12480 & 8.489 & 139.849 & 14.142 & 257.377 \\
		\hline
		$\Ft\Ct$ & 60 & 8.553 & 141.513 & 14.295 & 261.013 \\
		& 216 & 8.522 & 140.635 & 14.205 & 258.870 \\
		& 816 & 8.510 & 140.349 & 14.179 & 258.389 \\
		& 3168 & 8.507 & 140.257 & 14.300 & 262.429 \\
		& 12480 & 8.506 & 140.244 & 14.204 & 259.429 \\
	\end{tabular}
	\caption{Results for Cook's membrane example for $f=8$ and $f=16$ with $2\times2$, $4\times4$, $8\times8$, $16\times16$, and $32\times32$ grids. For all methods the number of coupling dofs, the vertical deflection at point $A$ and the $L^2$ norm of the displacement are presented.}
	\label{tab:cooks_membrane_result_f816}
\end{table}

\begin{table}[h!]
	\centering
	\begin{tabular}{cc|cc|cc}
		& coupl. dof & w & $\|U\|$& w & $\|U\|$ \\
		\hline
		&&$f=24$&&$f=32$&\\
		std & 32 & 13.305 & 212.422 & 16.348 & 265.766\\
		 & 112 & 15.579 & 269.195 & 18.574 & 329.713\\
		 & 416 & 17.052 & 315.428  & 20.053 & 385.690 \\
		 & 1600 & 17.700 & 337.958  & 20.778 & 415.150 \\
		 & 6272 & 17.948 & 346.933  & 21.093 & 427.745\\
		\hline
		$\Ft$ &  60 & 18.439 & 363.620 & 21.769 & 453.584\\
		 &  216 & 18.317 & 359.827 & 21.639 & 448.622\\
		 & 816 & 18.317 & 360.305 & 21.747 & 453.602\\
		 & 3168 & 18.346 & 361.587 & 21.612 & 448.416\\
		 & 12480 & 18.281 & 359.313 & 21.530 & 445.142\\
		\hline
		$\Ct$  & 60 & 17.961 & 345.705 & 21.217 & 429.143\\
		& 216 & 18.103 & 351.545 & 21.348 & 435.983\\
		& 816 & 18.157 & 354.038 & 21.402 & 439.090\\
		& 3168 & 18.188 & 355.592 & 21.453 & 441.644\\
		& 12480 & - & - & - & - \\
		\hline
		$\Ft\Ct$ & 60 & 18.424 & 362.626 & 21.734 & 451.528 \\
		& 216 & 18.292 & 358.944 & 21.596 & 446.903 \\
		& 816 & 18.584 & 370.337 & 21.769 & 455.071 \\
		& 3168 & 18.345 & 361.562 & 21.573 & 446.896\\
		& 12480 & 18.265 & 358.718  & 21.531 & 445.182\\
	\end{tabular}
	\caption{Results for Cook's membrane example for $f=24$ and $f=32$ with $2\times2$, $4\times4$, $8\times8$, $16\times16$, and $32\times32$ grids. For all methods the number of coupling dofs, the vertical deflection at point $A$ and the $L^2$ norm of the displacement are presented.}
	\label{tab:cooks_membrane_result_f2432}
\end{table}

\subsubsection{Thin Beam}
\label{subsubsec:thin_beam}
\begin{figure}[h!]
	\centering
	\includegraphics[width=0.75\linewidth]{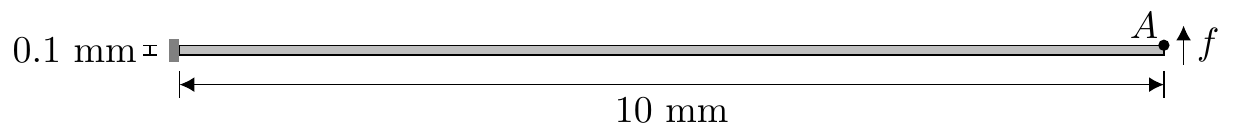}
	\caption{Geometry of thin beam example.}
	\label{fig:thin_beam_geo}
\end{figure}

For the thin beam example we follow \cite{R02,WBAR17}, where the beam is clamped at the left side and a point load at $A$ is applied. But instead of a point load we choose a shear force $f = \SI{1}{\newton\per\milli\meter\squared}$ on the right boundary such that the displacement is in the same magnitude, see Figure \ref{fig:thin_beam_geo}. This modification is needed as point-evaluation is not well defined in context of Sobolev spaces and traces, and cannot be set directly in terms of the TDNNS method. There, only the tangential trace is well defined for the displacement $\uv$ and the normal trace for the hybridization variable $\alphav$. The quantity of interest is again the vertical deflection at point $A$. Structured quadrilateral meshes with $10\times 1$, $20\times 2$, $40\times 4$, and $80\times 8$ grids are used. The material parameters are $\mu=\SI{6000}{\newton\per\milli\meter\squared}$ and $\lambda=\SI{24000}{\newton\per\milli\meter\squared}$ and the hyperelastic potential \eqref{eq:neohooke_log} is considered.

\begin{figure}[h!]
	\centering
	\includegraphics[width=0.5\linewidth]{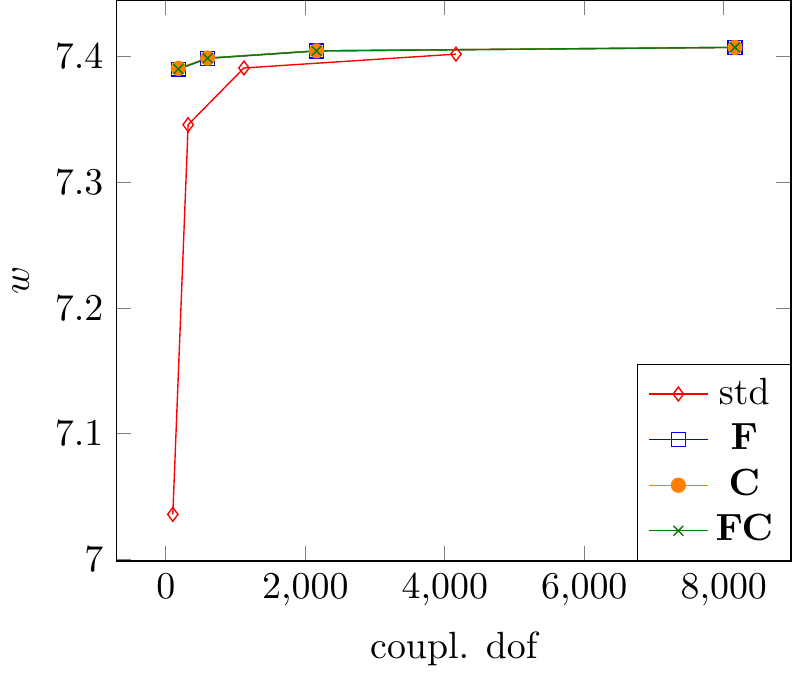}
	\caption{Vertical deflection of thin beam example at point $A$.}
	\label{fig:result_thin_beam}
\end{figure}

We observe that all three methods perform well, even for the coarsest grid. As discussed and shown in \cite{PS12} the linear TDNNS method is extremely robust for anisotropic elements. The standard Lagrangian method converges also, however, it needs a finer grid to produce accurate solutions, see  Table \ref{tab:result_thin_beam} and Figure \ref{fig:result_thin_beam} for the results.

\begin{table}[h!]
	\centering
	\begin{tabular}{c|ccc|ccc}
		&coupl. dof & w & $\|U\|$ &coupl. dof & w & $\|U\|$ \\
		\hline
		&ne=10 &&&ne=40&&\\
		std  & 100 & 7.036 & 3.992 & 320 & 7.346 & 4.273 \\
		$\Ft$   & 180 & 7.390 & 4.314	& 600 & 7.399 & 4.321 \\
		$\Ct$ &  180 & 7.391 & 4.314 & 600 & 7.399 & 4.321 \\
		$\Ft\Ct$   & 180 & 7.390 & 4.314 & 600 & 7.399 & 4.321 \\
		\hline
		&ne=160 &&&ne=640&&\\
		std & 1120 & 7.391 & 4.314 & 4160 & 7.402 & 4.324 \\
		$\Ft$ &  2160 & 7.405 & 4.326 & 8160 & 7.407 & 4.329 \\
		$\Ct$ & 2160 & 7.405 & 4.326 & 8160 & 7.407 & 4.329 \\
		$\Ft\Ct$ & 2160 & 7.405 & 4.326 &  8160 & 7.407 & 4.329 \\
	\end{tabular}
	\caption{Results for thin beam example. For all methods the number of coupling dofs, the vertical deflection at point $A$ and the $L^2$ norm of the displacement are given.}
	\label{tab:result_thin_beam}
\end{table}

\subsubsection{Inflation of a Hollow Spherical Ball}
\label{subsubsec:infl_sphere}
\begin{figure}[h!]
	\centering
	\includegraphics[width=0.4\linewidth]{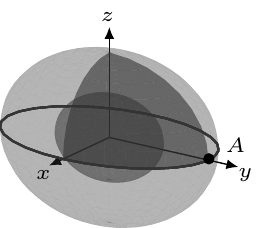}\hspace*{2cm}
	\includegraphics[width=0.2\linewidth]{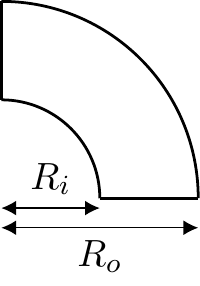}
	\caption{3D geometry of inflation of a hollow spherical ball example and 2D cross-section.}
	\label{fig:infl_sphere_geom}
\end{figure}

\begin{figure}[h!]
	\centering
	\includegraphics[width=0.35\linewidth]{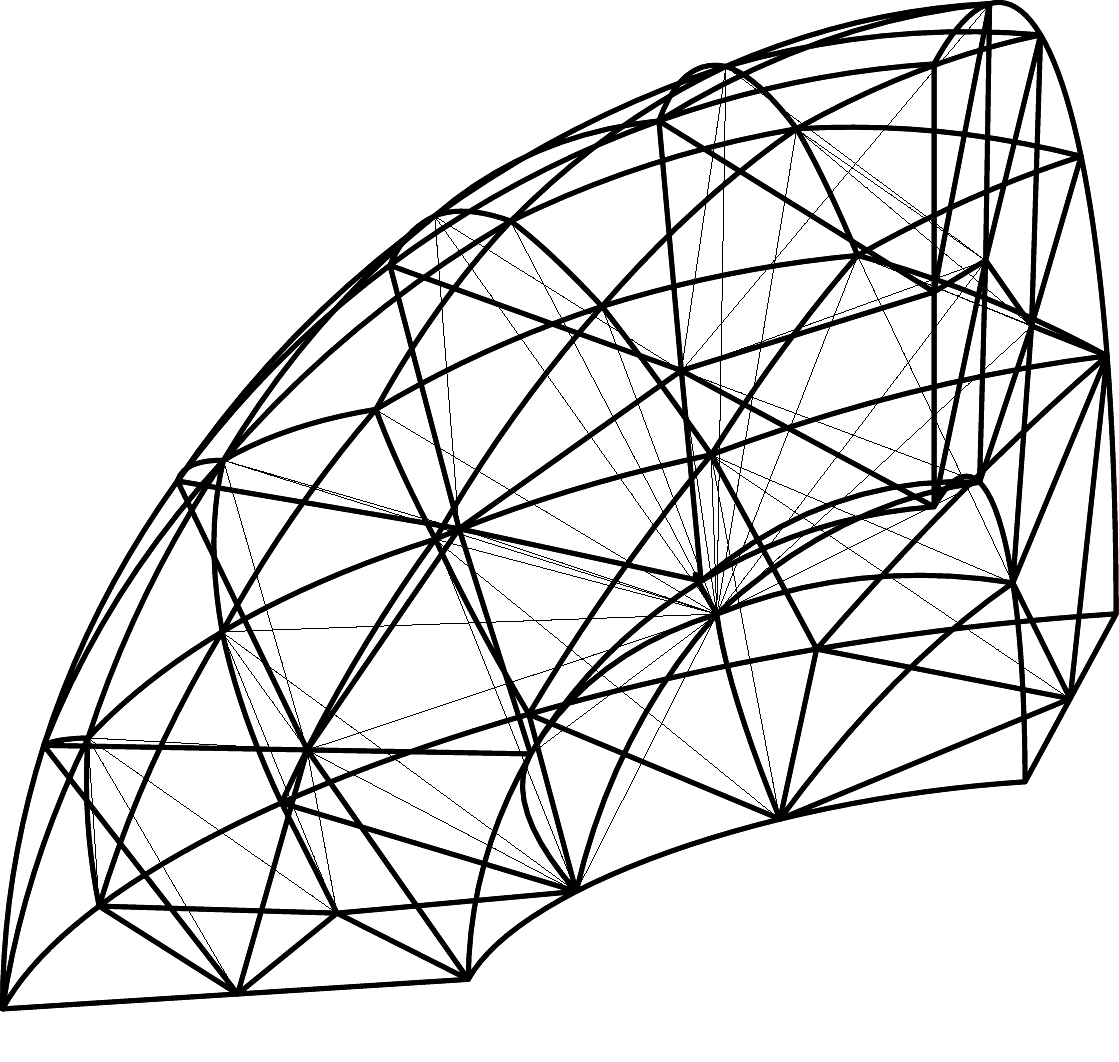}\hspace*{0.5cm}
	\includegraphics[width=0.45\linewidth]{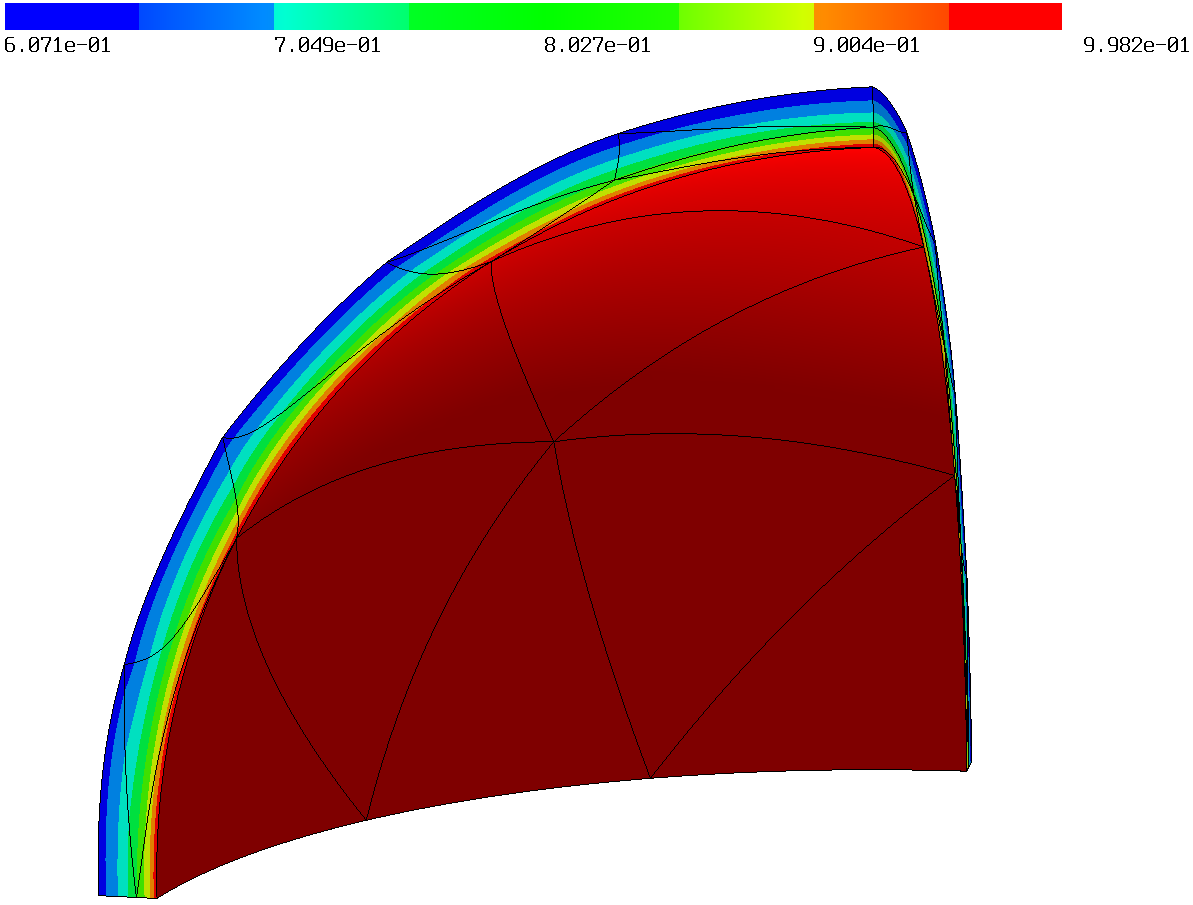}
	\caption{Curved mesh (left) and final configuration with 85 elements for inflation of a hollow spherical ball example.}
	\label{fig:infl_sphere_mesh}
\end{figure}

A hollow spherical ball is subjected to the boundary condition $\uv_{\text{in}}=(\gamma-1)\Xv$ on the inner boundary, whereas the outer boundary is left free. The inner and outer radius are given by $R_{\text{i}}=\SI{0.5}{\milli\meter}$ and $R_{\text{o}}=\SI{1}{\milli\meter}$, respectively. We adapted the benchmark in \cite{SY19} by setting the material parameters $\mu=\SI{1}{\newton\per\milli\meter\squared}$ and $\lambda = \SI{100}{\newton\per\milli\meter\squared}$ together with the hyperelastic potential \eqref{eq:neohooke_log} instead of using a complete incompressible material. This choice of $\lambda$, however, leads to a Poisson ratio $\nu =\frac{100}{202}\approx 0.495$ and thus, the example is already close to the nearly incompressible regime. The final configuration is reached for $\gamma=3$, starting from $\gamma=1$, the initial configuration.  Due to symmetry, only one eight of the ball is considered, see Figure \ref{fig:infl_sphere_geom}. Unstructured curved tetrahedral meshes are used as shown in Figure \ref{fig:infl_sphere_mesh}.

\begin{figure}[h!]
	\centering
	\includegraphics[width=0.5\linewidth]{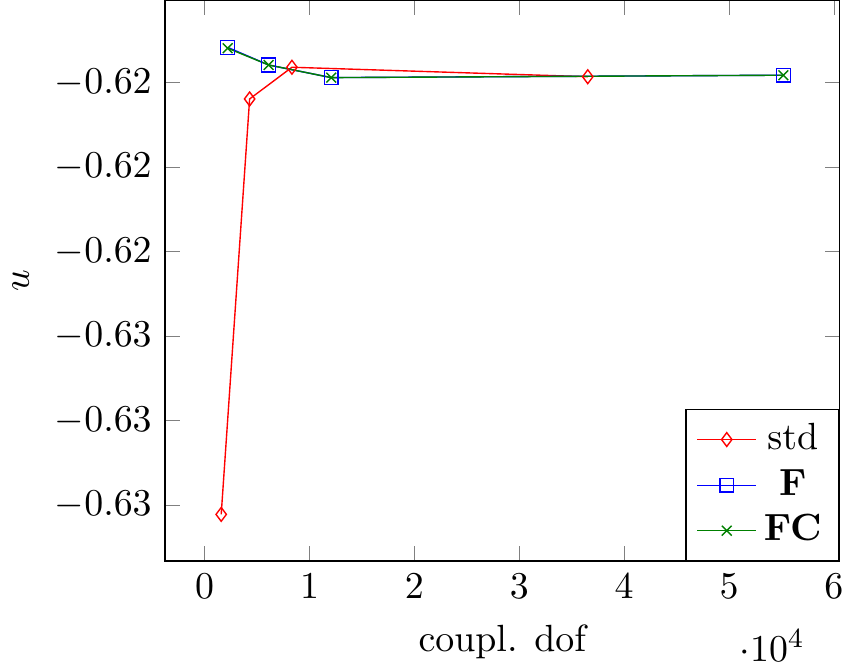}
	\caption{Radial deflection of inflation of a hollow spherical ball example at point $A$.}
	\label{fig:result_infl_sphere}
\end{figure}

We note that the standard method does not converge for polynomial order $2$, whereas the $\Ft$- and $\Ft\Ct$-based method do. Thus, only for the Lagrangian elements order $3$ is considered. The $\Ct$-based method, however, does not converge - even with the stabilization terms \eqref{eq:stab_u_alpha_n} and \eqref{eq:stab_C_Cu}. Although a higher polynomial order is used, the standard method does not perform as good as the $\Ft$- and $\Ft\Ct$-based method for the coarse grid, cf. Figure \ref{fig:result_infl_sphere}. The results are listed in Table \ref{tab:result_infl_sphere}.


\begin{table}[h!]
	\centering
	\input{infl_sphere.out}
	\caption{Results for inflation of a hollow spherical ball example. For all methods the number of elements, number of dofs and coupling dofs, the vertical deflection at point $A$ and the $L^2$ norm of the displacement are given.}
	\label{tab:result_infl_sphere}
\end{table}

\subsubsection{Cylindrical Shell}
\label{subsubsec:cyl_shell}

\begin{figure}[h!]
	\centering
	\includegraphics[width=0.35\linewidth]{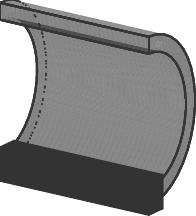}\hspace*{2cm}
	\includegraphics[width=0.22\linewidth]{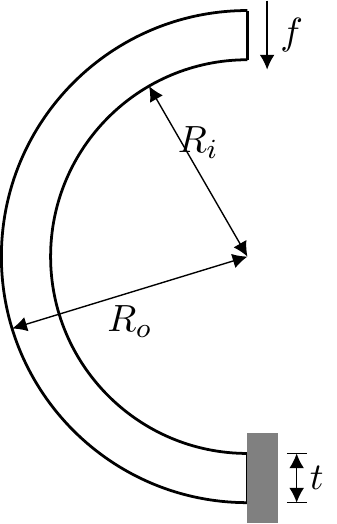}
	\caption{3D geometry of cylindrical shell example and 2D cross-section.}
	\label{fig:cyl_shell_geom}
\end{figure}

\begin{figure}[h!]
	\centering
	\includegraphics[width=0.33\linewidth]{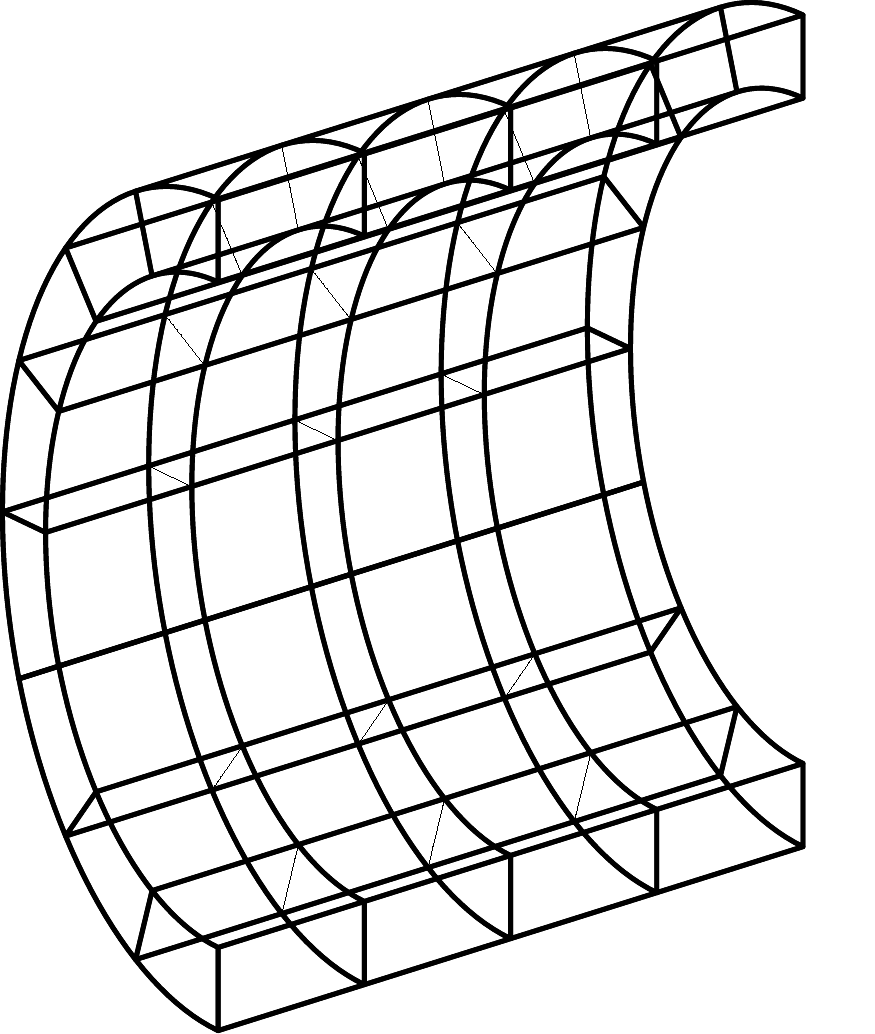}\hspace*{2cm}
	\includegraphics[width=0.39\linewidth]{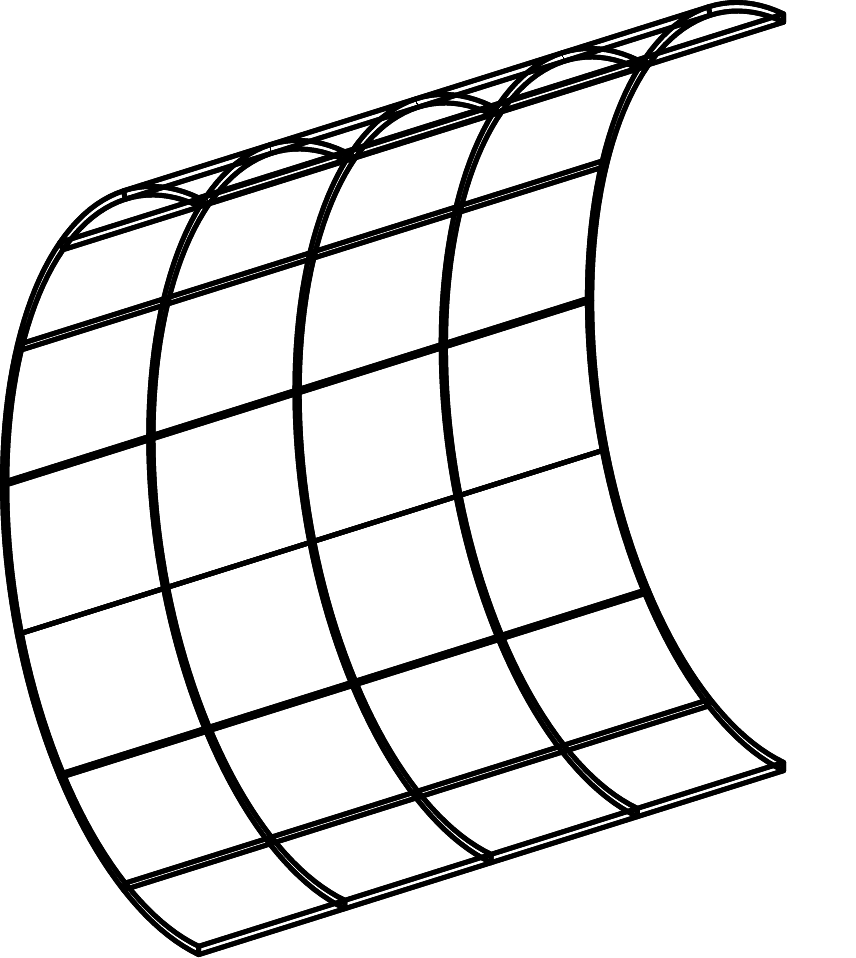}
	\caption{Curved meshes with $8\times 4\times1$ grid for cylindrical shell example for $t=2$ and $t=0.2$.}
	\label{fig:cyl_shell_mesh}
\end{figure}

\begin{figure}[h!]
	\centering
	\includegraphics[width=0.38\linewidth]{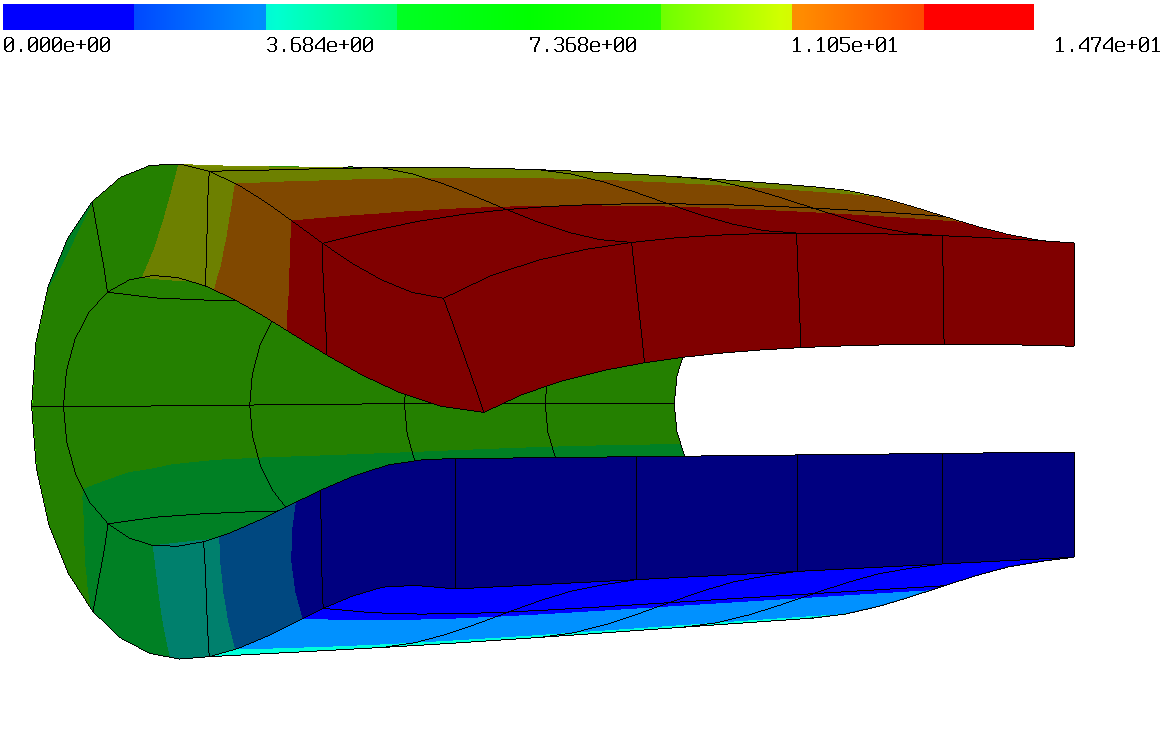}\hspace*{2cm}
	\includegraphics[width=0.38\linewidth]{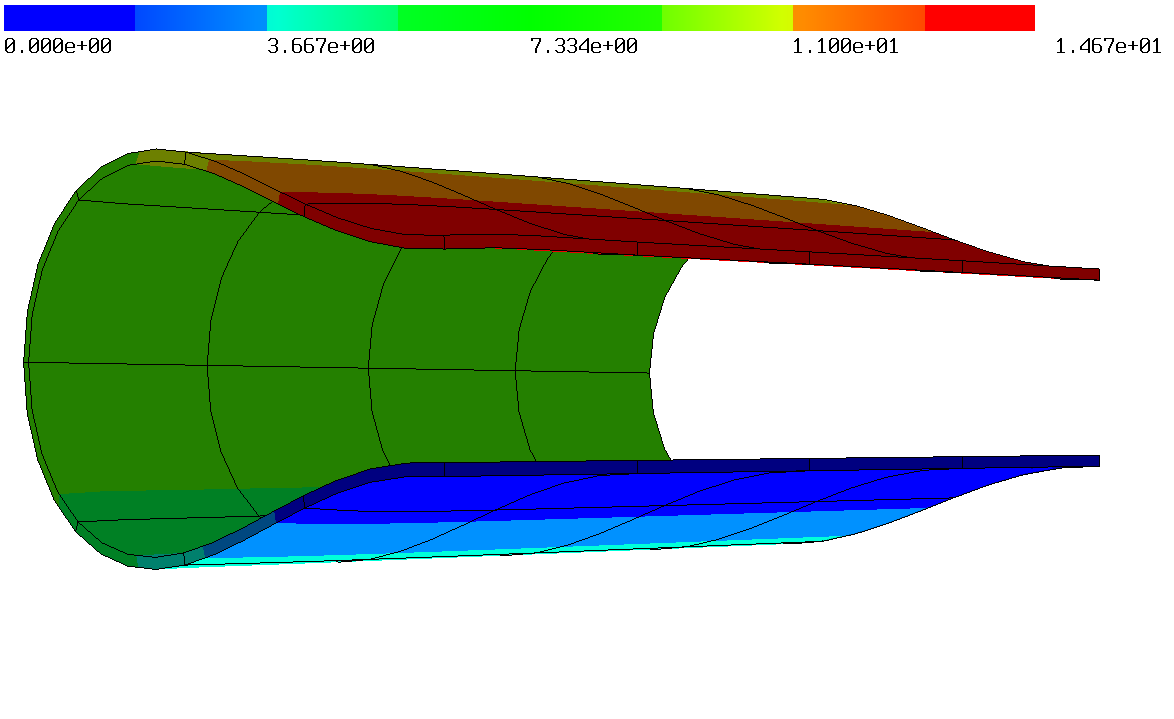}
	\caption{Final configuration of cylindrical shell example for $t=2$ and $t=0.2$ with $16\times 8\times1$ grid.}
	\label{fig:cyl_shell_final_def}
\end{figure}

\begin{figure}[h!]
	\centering
	\includegraphics[width=0.49\linewidth]{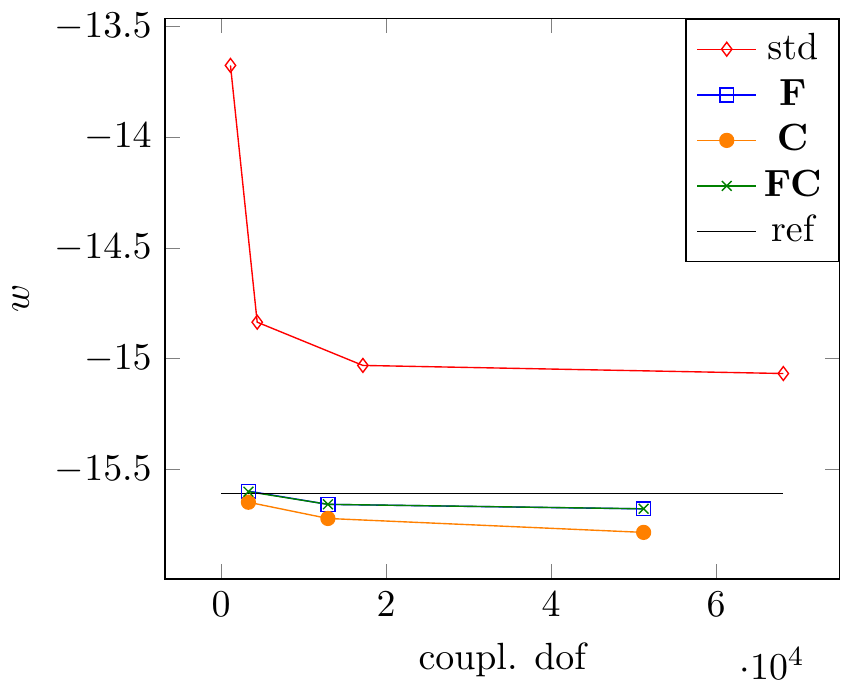}
	\includegraphics[width=0.47\linewidth]{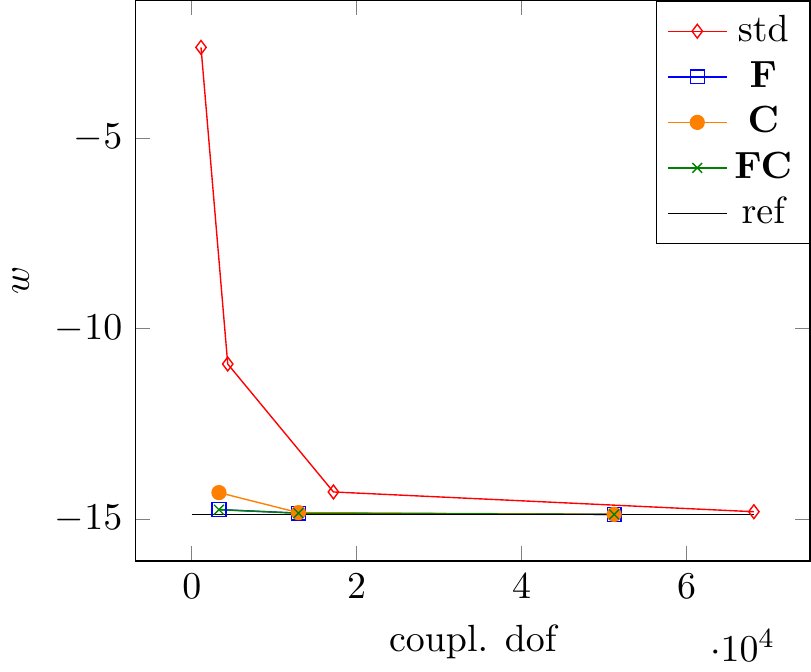}
	\caption{Vertical deflection of cylindrical shell example at point $A$.}
	\label{fig:result_cyl_shell}
\end{figure}

The benchmark presented in \cite{RWR98,R02} is adapted in terms of the force and boundary condition as line forces and traces are not well defined in terms of Sobolev spaces in three spatial dimensions. The same geometry and material parameters are considered. Namely, a quarter of a cylindrical structure with $r_i=\SI[parse-numbers = false]{9- t/2}{\milli\meter}$, $l=\SI{15}{\milli\meter}$, and thickness $t=\SI{2}{\milli\meter}$ or $\SI{0.2}{\milli\meter}$, see Figure \ref{fig:cyl_shell_geom}, and $\mu=\SI{6000}{\newton\per\milli\meter\squared}$ and $\lambda=\SI{24000}{\newton\per\milli\meter\squared}$ together with the hyperelastic potential \eqref{eq:neohooke_log}. The structure is clamped at the bottom area and an area shear force $f=\SI{240}{\newton\per\milli\meter\squared}$ ($f=\SI{2.7}{\newton\per\milli\meter\squared}$ for $t=\SI{0.2}{\milli\meter}$) is applied on the top. Structured hexahedra meshes with $8\times4\times1$, $16\times8\times1$, and $32\times16\times1$ grids are used, see Figure \ref{fig:cyl_shell_mesh}. For the standard method also a $64\times32\times1$ grid is used. The vertical deflection at point $A$ is depicted in Figure \ref{fig:result_cyl_shell}, the final deformations are shown in Figure \ref{fig:cyl_shell_final_def}, and the results can be found in Table \ref{tab:result_cyl_shell_t2}.  We observed a locking behaviour for the standard method, which becomes significant for the small thickness. All of the three presented methods give satisfying results already for the coarsest discretization. The reference values were computed with the standard method and degree $k=4$ on the finest grid, where locking is avoided due to the high polynomial degree. As only one layer in the thin direction is used the methods do not converge towards to reference solution for $t=2$. For the thin structure, $t=0.2$, the values match. As already mentioned the $\Ft$-based method may suffer from a larger number of Newton iterations to reach convergence, whereas the $\Ct$- and $\Ft\Ct$-based method have a better convergence behavior. In Table~\ref{tab:nit} the numbers of Newton iterations can be found, where twelve load steps were considered. Also the standard method needs more iterations than the $\Ct$- and $\Ft\Ct$-based variant.


\begin{table}[h!]
	\centering
	\begin{tabular}{c|c|cc|cc}
		& coupl. dof & w & $\|U\|$  & w & $\|U\|$\\
		\hline
		&&$t=2$&&$t=0.2$&\\ 
		\hline
		std  & 1125 & -13.676 & 211.958 & -2.613 & 14.350\\
		& 4365 & -14.835 & 231.071 & -10.927 & 58.065 \\
		& 17181 & -15.030 & 234.734 & -14.287 & 76.437\\
		& 68157 & -15.067 & 235.513 & -14.806 & 79.387\\
		\hline
		$\Ft$ &  3300 & -15.598 & 244.530  & -14.749 & 78.920\\
		& 12936 & -15.658 & 245.756 & -14.844 & 79.571 \\
		& 51216 & -15.678 & 246.149 & -14.876 & 79.825\\
		\hline
		$\Ct$ & 300 & -15.648 & 245.779  & -14.304 & 76.878\\
		& 12936 & -15.721 & 246.883 & -14.830 & 79.521\\
		& 51216 & -15.784 & 247.951 & -14.876 & 79.825 \\
		\hline
		$\Ft\Ct$ & 3300 & -15.600 & 244.586 & -14.753 & 78.946\\
		& 12936 & -15.658 & 245.757 & -14.844 & 79.576 \\
		& 51216 & -15.678 & 246.143  & -14.876 & 79.825 \\
	\end{tabular}
	\caption{Results for cylindrical shell example. For all methods the number of coupling dofs, the vertical deflection at point $A$ and the $L^2$ norm of the displacement are given.}
	\label{tab:result_cyl_shell_t2}
\end{table}

\begin{table}[h!]
	\centering
	\begin{tabular}{c|cccc}
		$t$&std & $\Ft$ & $\Ct$ & $\Ft\Ct$\\
		\hline
		$2$& 6 &5 &5&5\\
		$0.2$&7&7-8&4-5&4-5\\ 
	\end{tabular}
	\caption{Number of Newton iterations  for $32\times 16\times 1$ grid in cylindrical shell example.}
	\label{tab:nit}
\end{table}

\subsubsection{Thin Beam subjected to end moment}
\label{subsubsec:circle}
\begin{figure}[h!]
	\centering
	\includegraphics[width=0.48\linewidth]{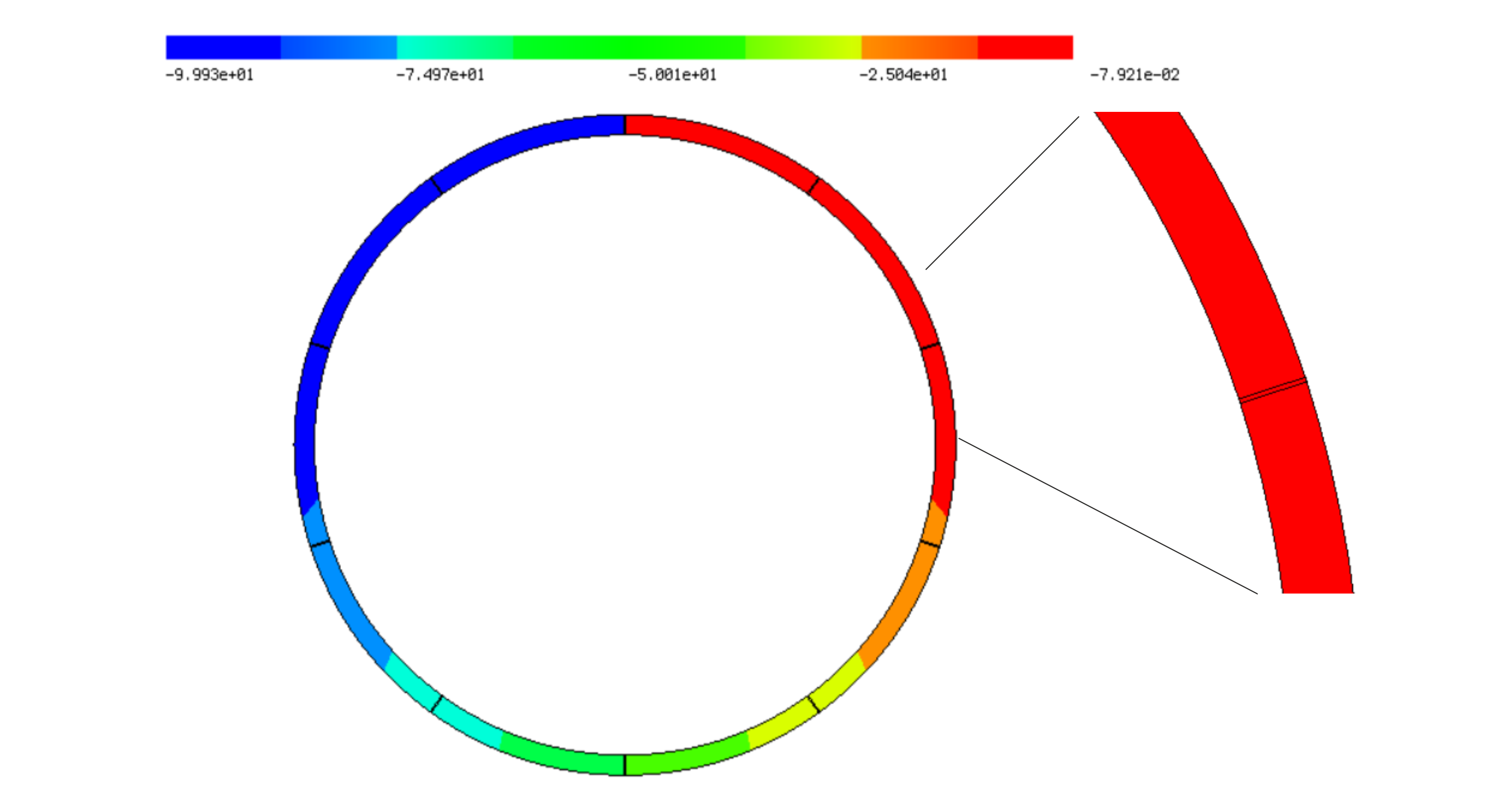}
	\includegraphics[width=0.48\linewidth]{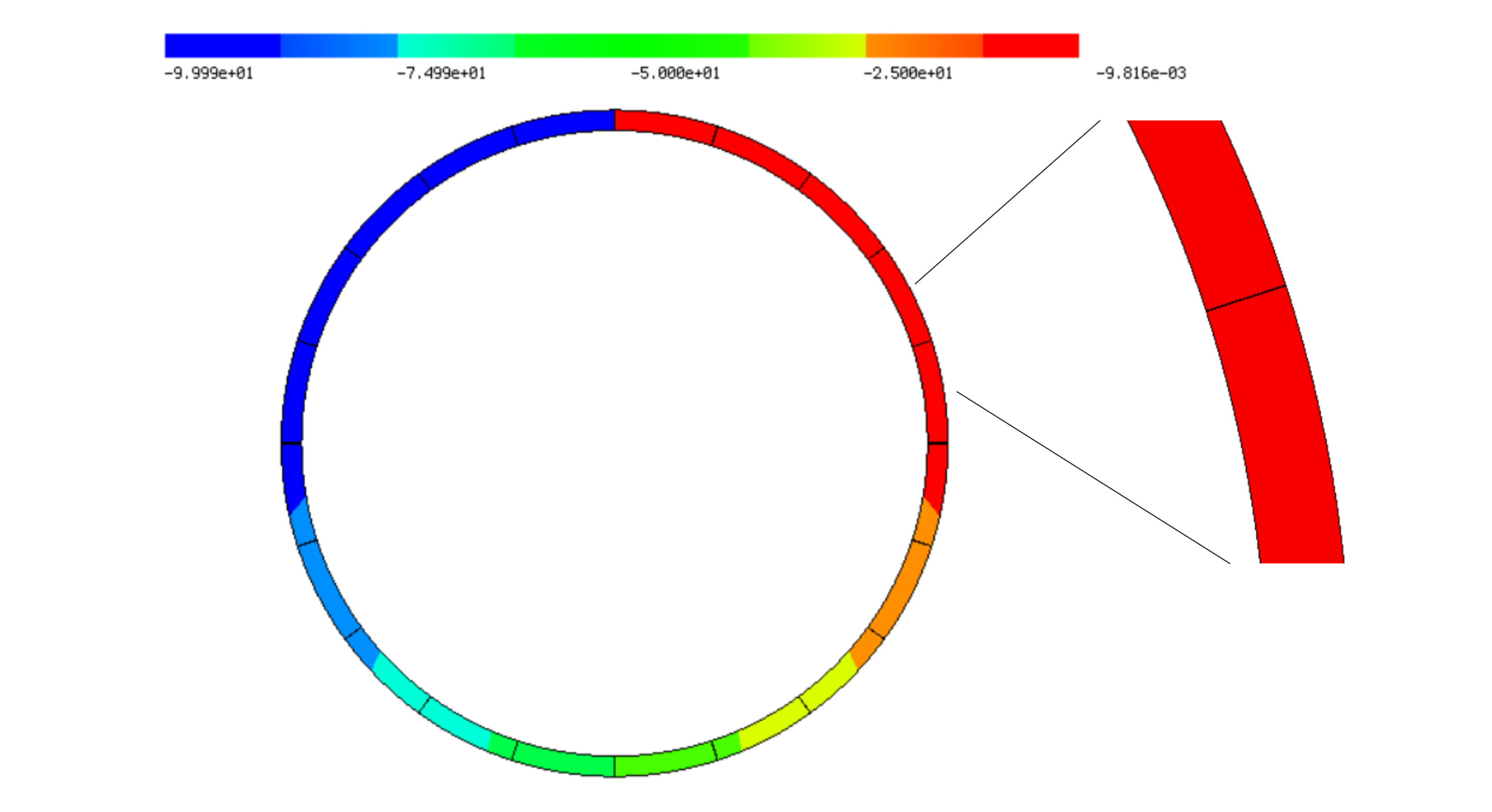}
	\caption{Displacement $u_x$ in deformed geometry of thin beam subjected to end moment,  $\Ft$-based method, including zoom into element interface. Left: $10$ elements, order $p=1$; right: $20$ elements, order $p=2$.}
	\label{fig:final_def_beam_mom}
\end{figure}

For this example we use a cantilever beam, clamped at the left-hand side, of length $L = \SI{100}{\milli\meter}$ and thickness $t = \SI{1}{\milli\meter}$. To allow for an anayltical solution, we choose material parameters $\mu=\SI{1e4}{\newton\per\meter\squared}$ and $\lambda=\SI{0}{\newton\per\meter\squared}$. An end moment is applied such that the beam cross-section rotates by $360^\circ$, thereby  forming a perfect circle in the theory of thin shells. Due to the enormous deformation the Updated Lagrangian scheme discussed in Section \ref{sec:updated_lagrangian} is used for all presented methods. Two different grids are used, where one element is used in thickness direction, and $10$ or $20$ elements in axial direction. The quantity of interest is given by the applied moment necessary to rotate the cross-section at the tip by $360^\circ$. In Figure \ref{fig:final_def_beam_mom} the final deformation is depicted and the results, as well as the absolute and relative error in the necessary moment as compared to the analytical value of $M_0 = \frac{4\pi \mu t^3}{12 L} = \SI{104.720}{\newton\meter}$ are given in Table \ref{tab:result_beam_mom}. 
The $\Ct$-based method did not converge for the coarsest grid and elements of order one, whereas the other two methods always converged. The overall accuracy of the methods is comparable. Concerning iteration counts, the $\Ct$- and $\Ft\Ct$-based methods converged very fast, needing 3-7 iterations in each of the 12 load steps. In comparison, for the standard nodal FEM and the $\Ft$-based methods, more than 20 iterations were necessary in some of the load steps.

\begin{table}[h!]
	\centering
	\begin{tabular}{l|ccc|ccc}
		& coupl. dof & moment &  rel. error & coupl. dof&  moment  & rel. error\\ \hline
	ne = 10	& order 1 &  & & order 2 & &\\
	std. method & 40 & 5538.381 & 51.888 & 100 & 231.149  & 1.207314\\
	$\Ft$ & 120 &  103.172  & -0.014781 & 180 & 104.685  & -0.000330 \\
	$\Ct$ & 120 & -- & -- & 180 & 104.744 & 0.000229 \\	
	$\Ft\Ct$ & 120 & 101.399  & -0.031715 & 180 & 104.679 & -0.000385 \\\hline
	ne = 20	& order 1 & & & order 2 & &\\	
	std & 80 & 1425.993  & 12.617 & 200 & 113.175  & 0.080743\\
	$\Ft$ & 240 & 103.908  & -0.007750 & 360 &  104.734  & 0.000138\\
	$\Ct$ & 240 & 104.668 & -0.000498 & 360 & 104.738 & 0.000177 \\
	$\Ft\Ct$ & 240 & 103.878  & -0.008038 & 360 & 104.734  & 0.000138 \\\hline
	\end{tabular}
	\caption{Results for circular bending of a thin beam. For all methods, the end moment necessary to obtain a rotation of $360^\circ$ is provided, as well as the relative error as compared to the analytical value of $M_0 = \frac{4\pi \mu t^3}{12 L}$.}
	\label{tab:result_beam_mom}
\end{table}

\section*{Acknowledgements}
\label{sec:acknowledgements}
The support by the Austrian Science Fund (FWF) projects W\,1245 and F\,65 is gratefully acknowledged.

\appendix

\section{Proof of consistency}
\label{app:consistency_C}
We prove consistency of problem \eqref{eq:lag_C1}. If we take the variation of \eqref{eq:lag_C2} in direction $\delta\uv$ and look at the left-hand side
\begin{align*}
\sum_{T\in\T}\left(\int_T\Ft(\uv)\Sigmat:\nabla\delta\uv\,d\Xv+\int_{\partial T}(\Ft(\uv)\Sigmat)_{\Nv\Nv}\delta\uv_{\Nv} - (\nabla\delta\uv\Sigmat)_{\Nv\Nv}(\uv-\alphav)_{\Nv}\,d\Sv\right) 
\end{align*}
we want to recover the strong form if we insert the true, smooth solution $\tilde{\uv}$. First, the hybridization variable $\alphav$ is the normal trace of $\tilde{\uv}$ and thus, $(\tilde{\uv}-\alphav)_{\Nv} = 0$ on $\partial T$. Integration by parts and reordering yields
\begin{align}
&\sum_{T\in\T}\left(\int_T-\text{div}(\Ft(\tilde{\uv})\Sigmat)\cdot\delta\uv\,d\Xv +\int_{\partial T}((\Ft(\tilde{\uv})\Sigmat)_{\Nv}\delta\uv-(\Ft(\tilde{\uv})\Sigmat)_{\Nv\Nv}\delta\uv_{\Nv}) \,d\Sv\right)\nonumber \\
=&\sum_{T\in\T}\left(\int_T-\text{div}(\Ft(\tilde{\uv})\Sigmat)\cdot\delta\uv\,d\Xv +\int_{\partial T}(\Ft(\tilde{\uv})\Sigmat)_{\Nv\Tv}\cdot\delta\uv_{\Tv} \,d\Sv\right)\nonumber\\
=&\sum_{T\in\T}\int_T-\text{div}(\Ft(\tilde{\uv})\Sigmat)\cdot\delta\uv\,d\Xv +\sum_{E\in\E}\int_{E}\llbracket(\Ft(\tilde{\uv})\Sigmat)_{\Nv_E\Tv_E}\rrbracket\cdot\delta\uv_{\Tv_E} \,d\Sv.
\end{align}

The first term states the element wise balance equation $-\text{div}(\Ft(\tilde{\uv})\Sigmat)|_T = \fv|_T$ and the second the continuity of the normal-tangential components of the first Piola--Kirchhoff stress tensor. Reordering hybridization terms in \eqref{eq:lag_C2} yields
\begin{align}
\sum_{E\in\E}\int_{E}\llbracket(\Ft(\tilde{\uv})\Sigmat)_{\Nv_E\Nv_E}\rrbracket\delta\alphav_{\Nv_E}\,d\Sv = 0
\end{align} 
forcing the normal-normal continuity and thus, the continuity of the normal component of the first Piola--Kirchhoff stress tensor. Hence, also the interface condition is fulfilled and the problem is consistent.

\section{Linearization}
\label{app:linear_C}
Computing the first variations and under the assumption of small deformations, i.e., $\uv=\mathcal{O}(\varepsilon)$, $\nabla\uv=\mathcal{O}(\varepsilon)$, $\alphav=\mathcal{O}(\varepsilon)$, $\Sigmat = \mathcal{O}(\varepsilon)$, and $\Ct(\uv)=2\epst(\uv)+\It+\mathcal{O}(\varepsilon^2)$, the variations of \eqref{eq:lag_C2} become 
\begin{subequations}
	\label{eq:var_prob_huwa_E_lin}
	\begin{alignat}{3}
	& \int_{\Omega}\frac{\partial\Psi(\Et)}{\partial \Et}:\delta \Et - \Sigmat:\delta\Et\,d\Xv &=&0\quad &\forall\delta\Et,\label{eq:var_prob_huwa_E_lin_a}\\
	&-\sum_{T\in\T}\left(\int_T(\Et-\epst(\uv)):\delta\Sigmat\,d\Xv+\int_{\partial T}(\delta\Sigmat)_{\Nv\Nv}(\uv-\alphav)_{\Nv}\,d\Sv\right) &=& 0 \quad &\forall \delta\Sigmat,\\
	&\sum_{T\in\T}\left(\int_T\Sigmat:\nabla\delta\uv\,d\Xv-\int_{\partial T}\Sigmat_{\Nv\Nv}\delta\uv_{\Nv} \,d\Sv\right) &=& \int_{\Omega}\fv\cdot\delta\uv\,d\Xv\quad &\forall\delta\uv,\\
	& \sum_{T\in\T}\int_{\partial T}\Sigmat_{\Nv\Nv}\delta\alphav_{\Nv}\,d\Sv &=& 0\quad &\forall \delta\alphav,
	\end{alignat}
\end{subequations}
where we implicitly defined $\Et:=\frac{1}{2}(\Ct-\It)$ and $\delta\Et:=\frac{1}{2}\delta\Ct$. Assuming a quadratic potential, i.e., $\frac{\partial\Psi}{\partial \Et} = \Mt\Et$, and eliminating $\Et$ by $\Sigmat$ with \eqref{eq:var_prob_huwa_E_lin_a} recovers the hybridized TDNNS method \eqref{eq:var_prob_lin_tdnns_hyb}. Thus, the linearized versions of \eqref{eq:lag_C2} and \eqref{eq:laghF1} coincide.

\bibliographystyle{acm}
\bibliography{cites}

\end{document}